\documentclass[11pt,reqno]{amsart}

\usepackage{amsthm,amsfonts,amsmath,amscd,amssymb,epsfig}
\usepackage[all]{xy}
\usepackage{enumerate,array}
\usepackage[latin9]{inputenc}
\usepackage{hyperref}
\usepackage{graphicx}
\usepackage{color}
\usepackage{amsfonts,amssymb}
\usepackage{mathrsfs}
\usepackage{a4wide}

\usepackage{color} 

\definecolor{green}{rgb}{0,0.5,0}
\definecolor{blue}{rgb}{0,0,1}

\hypersetup{colorlinks,
linkcolor=blue,
filecolor=green,
citecolor=green}

\advance\hoffset by -0.5 truecm

\theoremstyle{plain}

\newtheorem{neu}{}[section]
\newtheorem{Cor}[neu]{Corollary}
\newtheorem*{Cor*}{Corollary}
\newtheorem{Thm}[neu]{Theorem}
\newtheorem*{Thm*}{Theorem}
\newtheorem{Prop}[neu]{Proposition}
\newtheorem*{Prop*}{Proposition}
\theoremstyle{definition}
\newtheorem{Lemma}[neu]{Lemma}
\newtheorem*{Rmk*}{Remark}
\newtheorem{Rmk}[neu]{Remark}

\newtheorem*{Ex*}{Example}

\newtheorem*{Qu*}{Question}

\newtheorem{Def}[neu]{Definition}

\theoremstyle{remark}

\theoremstyle{definition}

\newcommand{\p}{\partial}
\newcommand{\om}{\omega}

\newcommand{\into}{\hookrightarrow}

\newcommand{\pf}{\longrightarrow}

\newcommand{\N}{{\mathbb{N}}}
\newcommand{\Z}{{\mathbb{Z}}}
\newcommand{\R}{{\mathbb{R}}}
\newcommand{\C}{{\mathbb{C}}}

\newcommand{\M}{\mathcal{M}}
\newcommand{\LLL}{\mathscr{L}}

\renewcommand{\H}{{\bf H}}

\renewcommand{\H}{\mathrm{H}}

\newcommand{\id}{\mathrm{id}}
\newcommand{\CF}{\mathrm{CF}}
\newcommand{\HF}{\mathrm{HF}}
\newcommand{\RFH}{\mathrm{RFH}}

\newcommand{\Crit}{{\rm Crit}}

\newcommand{\Ham}{\mathrm{Ham}}
\newcommand{\Supp}{\mathrm{Supp}}
\newcommand{\cl}{\mathrm{cl}}
\newcommand{\EE}{\mathcal{E}}

\newcommand{\FF}{\mathcal{F}}

\newcommand{\HH}{\mathcal{H}}
\newcommand{\GG}{\mathcal{G}}

\renewcommand{\AA}{\mathcal{A}}
\newcommand{\PP}{\mathcal{P}}

\newcommand{\FFF}{{\mathfrak F}}

\newcommand{\pp}{\mathfrak{p}}
\newcommand{\comment}[1]{}
\newcommand{\x}{\times}

\newcommand{\beq}{\begin{equation}}
\newcommand{\beqn}{\begin{equation}\nonumber}
\newcommand{\eeq}{\end{equation}}

\newcommand{\bea}{\begin{equation}\begin{aligned}}
\newcommand{\bean}{\begin{equation}\begin{aligned}\nonumber}
\newcommand{\eea}{\end{aligned}\end{equation}}

\numberwithin{equation}{section}
\begin{document}
\title [Rabinowitz Floer homology and coisotropic intersections] {Generalized Rabinowitz Floer homology and coisotropic intersections}
\author{Jungsoo Kang}

\address{Department of Mathematics and Research Institute of Mathematics\\
     Seoul National University}
\email{hoho159@snu.ac.kr}

\begin{abstract}
In this paper, we extend Rabinowitz Floer homology theory which has been established and extensively studied for hypersurfaces to coisotropic submanifolds of higher codimension. With this generalized version of Rabinowitz Floer homology theory, we explore the coisotropic intersection problem which interpolates between the Lagrangian intersection problem and the closed orbit problem. To be specific, we study the existence of leafwise intersection points on contact coisotropic submanifolds and the displaceability of stable coisotropic submanifolds.
\end{abstract}
\maketitle

\setcounter{tocdepth}{1}
\tableofcontents

\section{Introduction and main results}\label{introduction}
The coisotropic intersection problem was first studied in depth by Ginzburg \cite{Gi1}, and have been explored by many authors, see subsection \ref{History and related results}. Rabinowitz Floer homology theory which was developed by Cieliebak-Frauenfelder \cite{CF} using the action functional introduced by Rabinowitz \cite{Ra} is one of the effective methods to study the coisotropic intersection problem for hypersurfaces. We extend this theory to coisotropic submanifolds of arbitrary codimension. In this paper, we explore the existence problem of leafwise intersection points on (restricted) contact coisotropic submanifolds and the displaceability of stable coisotropic submanifolds. Furthermore, we define the generalized Rabinowitz Floer homology and compute it in the easiest case.\\[-2ex]

Throughout this paper, we deal with a symplectically aspherical symplectic $2n$-dimensional manifold $(M,\om)$ which is either closed or convex at infinity. We call a symplectic manifold $(M,\om)$ {\em convex at infinity} if $(M,\om)$ is symplectomorphic to the positive part of the symplectization of a closed contact manifold at infinity; $(M,\om)$ is called {\em symplectically aspherical} if one has the equality $\om|_{\pi_2(M)}=0$. Let $\Sigma$ be a closed submanifold in $(M,\om)$ of codimension $k$. We define the symplectic orthogonal space of $\Sigma$ by $T\Sigma^\om=\ker\om|_\Sigma$. Then $\Sigma$ is said to be {\em coisotropic} if $T_x\Sigma^\om$ is a $k$-dimensional subspace of $T_x\Sigma$ for every $x\in\Sigma$. We note that $k\in[0,n]$. The notions of stable, contact, and restricted contact type for coisotropic submanifolds were introduced by Bolle \cite{Bo1,Bo2}.

\begin{Def}
A coisotropic submanifold $\Sigma$ of codimension $k$ in $(M,\om)$ is called {\em stable} if there exist 1-forms $\alpha_1,\dots,\alpha_k$ on $\Sigma$ which satisfy
\begin{enumerate}
\item $\ker d\alpha_i\supset\ker\om_\Sigma$ for $i=1,\dots,k$;
\item $\alpha_1\wedge\cdots\wedge\alpha_k\wedge\om_\Sigma^{n-k}\ne0$.
\end{enumerate}
We say that $\Sigma$ is of {\em contact type} if $\alpha_1,\dots,\alpha_k$ are primitives of $\om_\Sigma$. If there are 1-forms $\lambda_1,\dots,\lambda_k$ on $M$ such that $d\lambda_i=\om$ and $\lambda_i|_\Sigma=\alpha_i$ for all $i=1,\dots,k$,  $\Sigma$ is said to be of {\em restricted} contact type.
\end{Def}

\noindent\textbf{Assumptions on manifolds}.
\begin{itemize}
\item A coisotropic submanifold $\Sigma$ in $(M,\om)$ is closed and of stable or contact or restricted contact type.
\item If a coisotropic submanifold $\Sigma$ is of stable or contact type, we assume that $(M,\om)$ is symplectically aspherical and either closed or convex at infinity.
\item If $\Sigma$ is a restricted contact coisotropic submanifold, $(M,\om)$ is automatically symplectically aspherical but never closed; so, if this is the case $(M,\om)$ is only assumed to be convex at infinity.\\[-1.5ex]
\end{itemize}

\noindent\textbf{Convention and Notations.}
\begin{itemize}
\item The {\em Hamiltonian vector field} $X_F$ associated to a Hamiltonian function $F\in C^\infty(S^1\x M)$ is defined implicitly by $i_{X_F}\om=dF$.
\item The flow of $X_F$ is denoted by $\phi_F^t$. The time one map of the flow $\phi_F=\phi_F^1$ is called a {\em Hamiltonian diffeomorphism}.
\item We denote by $\Ham_c(M,\om)$ the group of Hamiltonian diffeomorphisms generated by compactly supported Hamiltonian functions.
\item We denote by $||F||$ resp. $||\phi||$ the Hofer-norm of $F\in C_c^\infty(S^1\x M)$ resp. $\phi\in\Ham_c(M,\om)$ and they will be defined in Section 4.
\end{itemize}

In Section \ref{coisotropic}, we define the Reeb vector fields $R_1,\dots,R_k$ on stable coisotropic submanifolds of codimension $k$. The following two equations \eqref{Bolle's eq} and \eqref{equation} play crucial roles in the generalized Rabinowitz Floer theory and the coisotropic intersection problem. For $v\in C^\infty (S^1,\Sigma)$,
\beq\tag{Bo}\label{Bolle's eq}
\partial_t v(t)=\sum_{i=1}^k \eta_i R_i(v(t)),\qquad t\in S^1=\R/\Z,\,\, \eta_i\in\R.
\eeq\\[-1ex]
We note that constant loops $v$ in $\Sigma$ are trivial solutions of \eqref{Bolle's eq} with $\eta_i=0$ for all $i=1,\cdots,k$. Solutions of  \eqref{Bolle's eq} can be viewed as the generalized version of Reeb orbits; when $k=1$, solutions are nothing but Reeb orbits on contact manifolds. This equation \eqref{Bolle's eq} was first studied by Bolle. We show in Lemma \ref{Lem:eta's are same} that \eqref{Bolle's eq} can be reduced to the following equation \eqref{equation} in the restricted contact case.

\beq\tag{Eq}\label{equation}
\partial_t v(t)=\eta\sum_{i=1}^k  R_i(v(t)), \qquad t\in S^1=\R/\Z,\,\, \eta\in\R.
\eeq\\[-1ex]
It is noteworthy that $\eta$ becomes the period of $v_\eta(t):=v(t/\eta)$, $t\in\R/\eta\Z$ solving
\beq\label{gen. Reeb eq.}
\p_t v_\eta(t)=\sum_{i=1}^k R_i(v_\eta(t))
\eeq
when $(v(t),\eta)$ solves \eqref{equation}, see Lemma \ref{Lem:eta's are same} together with \eqref{eq:eta is period}.\\[-1ex]

\subsection{Leafwise intersections}
Let $(M,\om)$ be a $2n$-dimensional symplectic manifold and $\Sigma$ be a closed coisotropic submanifold of codimension $k$. Then the symplectic structure $\om$ determines the symplectic orthogonal bundle $T\Sigma^\om\subset T\Sigma$ as follows:
$$
T\Sigma^\om:=\big\{(x,\xi)\in T\Sigma\,\big|\,\om_x(\xi,\zeta)=0 \textrm{ for all } \zeta\in T_x\Sigma\big\}.
$$
Since $\om$ is closed, $T\Sigma^\om$ is integrable, thus $\Sigma$ is foliated by leaves of the characteristic foliation. We denote by $L_x$ the leaf through $x$. A point $x\in\Sigma$ is called a {\em leafwise (coisotropic) intersection point} of $\phi\in\Ham_c(M,\om)$ if $\phi(x)\in L_x$. In the extremal case $k=n$, a coisotropic submanifold is foliated by only one leaf and is a Lagrangian submanifold. Thus a Lagrangian intersection point coincides with a leafwise intersection point in the Lagrangian case. We note that if there are Poisson-commute Hamiltonian functions $G_1,\dots,G_k\in C^\infty(M)$ (see Definition \ref{def:Poisson-commute}) which have $0$ as a regular value, then $\Sigma=\bigcap_{i=1}^k G_i^{-1}(0)$ is a coisotropic submanifold in $M$ and the leaf $L_x$, $x\in\Sigma$ can be written by
$$
L_x=\{\phi_{G_1}^{t_1}\circ\phi_{G_2}^{t_2}\circ\cdots\circ\phi_{G_k}^{t_k}(x)\,|\,t_1,\dots t_k\in\mathbb{R}\}.
$$

Albers-Frauenfelder \cite{AF1} showed that the perturbed Rabinowitz action functional and Rabinowitz Floer homology are well suited to the leafwise intersection problem for restricted contact hypersurfaces. One tough assumption that they required is the {\em separating condition} for hypersurfaces. That is, a hypersurface $\Sigma$ separates $M$ into two connected components of which one is compact. Without the separating condition, it is impossible to find a defining Hamiltonian function $G\in C^\infty(M)$ of $\Sigma$ such that $G^{-1}(0)=\Sigma$ in general.

In the same vein, when $\Sigma$ is a coisotropic submanifold of codimension $k$, we need the following generalized separating condition to generalize Rabinowitz Floer homology theory: There exist Hamiltonian functions $G_1,\cdots, G_k\in C^\infty(M)$ of $\Sigma$ such that $\bigcap_{i=1}^kG_i^{-1}(0)=\Sigma$ (to be called global coordinates, see Definition \ref{Def:global Bolle's coordinates}). But, at present, there are a few examples satisfying the (generalized) separating condition and thus we prove some results in the present paper without the separating condition. Unfortunately, however, it is still indispensable to define Rabinowitz Floer homology.

\begin{Def}
We denote by $\wp(\Sigma)>0$ the {\em minimal symplectic area} of all solutions of \eqref{Bolle's eq} contractible in $M$. To be more exact,
$$
\wp(\Sigma):=\inf\big\{|\Omega(v)>0|\,\big|\,v\,\,\textrm{solving \eqref{Bolle's eq} and contractible in $M$},\,v(t)\in\Sigma,\,t\in S^1\big\}.
$$
Here $\Omega$ stands for the symplectic area functional, i.e.  $\Omega(v)=\int_{D^2}\bar v^*\om$ where $\bar v\in C^\infty(D^2,M)$ is a filling disk of $v$, i.e. $\bar v|_{\p D^2}(t)=v(t)$ for $t\in S^1$. The symplectic asphericity condition guarantees that the value of $\Omega(v)$ is independent of the choice of a filling disk. When $\Sigma$ is of restricted contact type, then \eqref{Bolle's eq} is reduced to \eqref{equation} and $\wp(\Sigma)$ becomes the {\em minimal period} for solutions of \eqref{gen. Reeb eq.}. If there are no solutions of  \eqref{Bolle's eq}, we set $\wp(\Sigma)=\infty$ by convention.
\end{Def}

\noindent\textbf{Theorem A.}\textit{ Let $\Sigma$ be a closed restricted contact coisotropic submanifold in a symplectic manifold $(M,\om)$ being convex at infinity. If $||\phi||<\wp(\Sigma)$, there exists a leafwise intersection point for $\phi\in\Ham_c(M,\om)$ .}

\begin{Rmk}
G\"urel \cite{Gu} also proved Theorem A using a different method. We cannot entirely drop the restricted contact condition in Theorem A, see \cite[Example 7.2]{Gi1} and \cite[Remark 1.4]{Gu}. On the existence of infinitely many leafwise intersection points, we refer to \cite{AF2,AF3,Ka2,Ka3}.
\end{Rmk}

If a contact coisotropic submanifold is not of restricted contact type, our ambient symplectic manifold need not be exact and can be closed. So we have more examples. Even if a coisotropic submanifold $\Sigma$ is of contact type, we still can find a leafwise intersection point for restricted perturbations. Before stating the result, we introduce some notation. We will meet these again in Section \ref{coisotropic} and Section \ref{unrestricted}. We set
$$
U_r:=\bigr\{(q,\pp)=(q,p_1,\dots,p_k)\in\Sigma\x\R^k\,\bigr|\,|p_i|<r,\,\, \textrm{for all } i=1,\dots,k\bigr\},
$$
$$
\delta_0:=\max\bigr\{r\in\R\,\bigr|\,\textrm{there exists a symplectic embedding } \,\psi:U_r\into M\bigr\}.
$$
Let $\psi_0:U_{\delta_0}\into M$ be a symplectic embedding. Throughout this paper, we tacitly identify $U_{\delta}$ with $\psi_0(U_{\delta})$ for all $0<\delta\leq\delta_0$. For a time dependent Hamiltonian function $F\in C_c^\infty(S^1\x M)$, we define the support of the Hamiltonian vector field $X_F$ as
$$
\Supp X_F:=\big\{x\in M\,\big|\,X_F(t,x)\ne0 \textrm{ for some }t\in S^1\big\}.\\[0.5ex]
$$
We call a Hamiltonian function $F\in C_c^\infty(S^1\x M)$ {\em admissible} if $F$ is constant outside of $U_{\delta_0}$, i.e. $\Supp X_F\subsetneq U_{\delta_0}$. We denote by $\FFF$ the set of all admissible Hamiltonian functions:
$$
\FFF:=\big\{F\in C_c^\infty(S^1\x M)\,|\,\Supp X_F\subsetneq U_{\delta_0}\big\}.
$$

\noindent\textbf{Theorem B.}\textit{ Let $\Sigma$ be a closed contact coisotropic submanifold in a symplectically aspherical symplectic manifold $(M,\om)$ which is either closed or convex at infinity. Then $\phi_F$ for $F\in\FFF$ has a leafwise intersection point provided  $||F||<\wp(\Sigma)$.}\\[-1.5ex]

\noindent\textbf{Corollary B.}\textit{ The standard torus $(S^1)^n$ embedded in $(\C^*)^n$ with the standard symplectic structure $\om_{\mathrm{std}}$ has a self-intersection point for every $\phi\in\Ham_c((\C^*)^n,\om_{\mathrm{std}})$.}

\begin{Rmk}
$(S^1)^n\subset(\C^*)^n,\om_{\mathrm{std}})$ is certainly of restricted contact type. But the technique used to proved Corollary B is similar to the proof of Theorem B. Even though the symplectic manifold $(\C^*)^n,\om_{\mathrm{std}})$ is {\em not} convex at infinity, we can show that  gradient flow lines of a parametrized perturbed Rabinowitz action functional with some asymptotic condition do not escape to infinity. Then the proof of Theorem A guarantees the existence of a self-intersection point of $(S^1)^n$ for $\phi\in\Ham_c((\C^*)^n,\om_{\mathrm{std}})$ with $||\phi||<\wp((S^1)^n)=\infty$.
\end{Rmk}

\noindent\textbf{Theorem C.}\textit{ Let $\Sigma$ be a restricted contact coisotropic submanifold with global coordinates in a symplectic manifold $(M,\om)$. Then the number of leafwise intersection points for a generic  $\phi\in\Ham_c(M,\om)$ with $||\phi||<\wp(\Sigma)$ is bounded below by the sum of $\Z/2$-betti numbers of $\Sigma$.}\\[-1.5ex]

Theorem A and Theorem C were proved by Albers-Frauenfelder \cite{AF1} for restricted contact hypersurfaces with the separating condition. One thing to remark is that in Theorem A we succeed in removing the separating condition. But proving Theorem C we make use of Rabinowitz Floer homology. As we mentioned, the separating condition is indispensable to define Rabinowitz Floer homology and hence in Theorem C.\\[-1.5ex]

\subsection{Displacement energy}
A submanifold $\Sigma$ in a symplectic manifold $(M,\om)$ is said to be {\em displaceable} if there exists a Hamiltonian diffeomorphism $\phi\in\Ham_c(M,\om)$ such that $\phi(\Sigma)\cap\Sigma=\emptyset$. The {\em displacement energy} of $\Sigma$ in $M$ is defined by
$$
e(\Sigma):=\inf\big\{||F||\,\big|\,F\in C^\infty_c(S^1\x M),\,\,\phi_F(\Sigma)\cap\Sigma=\emptyset\big\}.
$$
We set $e(\Sigma)=\infty$ for the infimum of the empty set; that is, the displacement energy of a nondisplaceable submanifold is infinity.\\[-1.5ex]

\noindent\textbf{Theorem D.}
\textit{Let $\Sigma$ be a displaceable closed stable coisotropic submanifold in $(M,\om)$ which is closed (or convex at infinity) and symplectically aspherical. Then there exists a solution $v\in C^\infty(S^1,\Sigma)$ of \eqref{Bolle's eq} contractible in $M$, such that}
\beq\label{eq:Theorem D}
0<|\Omega(v)|\leq e(\Sigma).\\[1ex]
\eeq
\begin{Rmk}
The estimation \eqref{eq:Theorem D} is sharp: The unit sphere $S^{2n-1}\subset\R^{2n}$ has $e(S^{2n-1})=\pi=\Omega(v)$ where $v$ is a Reeb orbit of the standard contact structure on $S^{2n-1}$. For displaceable closed restricted contact coisotropic submanifolds, Theorem D was proved by Ginzburg \cite{Gi1}. A similar result was also proved by Cieliebak-Frauenfelder-Paternain \cite{CFP} for stable hypersurfaces using Rabinowitz Floer theory. Making use of their proof, we slightly improve their theorem.
\end{Rmk}

Let an ambient symplectic manifold $M$ be exact with a symplectic structure $\om=d\lambda$. Then $\lambda|_\FF=\lambda_\FF$, the restriction of $\lambda$ to the characteristic foliation $\FF$ of $\Sigma$ is closed. The cohomology class $[\lambda_\FF]\in\H^1_\mathrm{dR}(\FF)$ in the foliated de Rham cohomology (see \cite{MoSc}) is called the {\em coisotropic Liouville class}. In the case that the dimension of $\Sigma$ is half the dimension of $M$, it coincides with the ordinary Liouville class defined on Lagrangian submanifolds, see \cite{Po}. This generalized version of the Liouville class was considered by Ginzburg \cite{Gi1} and he also deduced Corollary D. For a given solution $v\in C^\infty(S^1,\Sigma)$ of \eqref{Bolle's eq} contractible in $M$, choose any filling disk $\bar v:D^2\to M$, i.e. $\bar v|_{\p D^2}(t)=v(t)$. In the sense of the following formula, we refer to $(\lambda,v)$ as the {\em symplectic area of $v$}.
$$
(\lambda,v)=\int_{S^1}v^*\lambda_\FF=\int_{D^2}\bar v^*\om.
$$\\[0.3ex]
Accordingly, to show that $[\lambda_\FF]\ne0$ it is equivalent to find a loop tangent to the foliation which has nonzero symplectic area. Therefore Corollary D below is an immediate consequence of Theorem D.\\[-1.5ex]

\noindent\textbf{Corollary D.}
\textit{A closed stable coisotropic submanifold $\Sigma$ in an exact symplectic manifold $(M,d\lambda)$ being convex at infinity is non-displaceable provided that  $[\lambda_\FF]=0$.}

\subsection{Rabinowitz Floer homology}
There are two types of Rabinowitz action functionals. We denote by $\AA^\HH$ the (unperturbed) Rabinowitz action functional  which is generically Morse-Bott. Here $\HH\in C^\infty(M,\R^k)$ and $\HH^{-1}(0)$ is a coisotropic submanifold that we would like to investigate. The chain complex for Floer homology of $\AA^\HH$ is generated by critical points of an auxiliary Morse function on the solution space of \eqref{Bolle's eq} and the boundary map is defined by counting gradient flow lines of the Morse function with cascades of $\AA^\HH$ (based on Frauenfelder's Morse-Bott homology \cite{F}). On the other hand, the perturbed Rabinowitz action functional is generically Morse and denoted by $\AA^\HH_F$ where $F\in C^\infty(S^1\x M,\R)$. The chain complex for Floer homology of $\AA^\HH$ is generated by leafwise intersection points and the boundary map is defined by counting gradient flow lines of $\AA^\HH_F$. Here cascades resp. gradient flow lines of $\AA^\HH$ resp. $\AA^\HH_F$ are solutions of a nonlinear elliptic PDE.

One of the power of Floer homology is the invariance property. Two Floer homologies obtained by $\AA^\HH$ and $\AA^\HH_F$ are isomorphic due to the standard continuation argument in Floer theory, see Section \ref{Rabinowitz Floer Homology}. Thus we name Rabinowitz Floer homology for both and denote by
$$
\RFH(\Sigma,M):=\HF(\AA^\HH)\cong\HF(\AA^\HH_F).
$$
We should mention that $\RFH(\Sigma,M)$ does not depend on the choice of  $\HH\in C^\infty(M,\R^k)$ the defining Hamiltonian tuple for $\Sigma$ (up to canonical isomorphism).

\begin{Rmk}\label{rmk:RFH for contact and stable}
We again emphasize that the separating condition is necessary for Rabinowitz Floer homology. Though we only deal with restricted contact coisotropic submanifolds, it is possible to define $\HF(\AA^\HH)$ in the stable case or $\HF(\AA_F^\HH)$ with $F\in\FFF$ in the contact case. The assertions (i) and (ii) in Theorem E continue to hold for contact coisotropic submanifolds if we restrict the class of perturbations to $\FFF$ and (iii) holds true for stable coisotropic submanifolds.
\end{Rmk}
The following theorem is an immediate consequence of the construction and invariance property of Rabinowitz Floer homology. \\[-1.5ex]

\noindent\textbf{Theorem E.} \textit{Let $\Sigma$ be a closed restricted contact coisotropic submanifold with global coordinates in a symplectic manifold $(M,\om)$ being convex at infinity. }
\begin{enumerate}
\item \textit{If Rabinowitz Floer homology does not vanish, there exists a leafwise intersection point for every $\phi\in\Ham_c(M,\om)$. In particular, if $\Sigma$ is displaceable in $M$, Rabinowitz Floer homology vanishes.}

\item \textit{There exists a nontrivial solution of \eqref{Bolle's eq} contractible in $M$, provided that $\Sigma$ is displaceable in $M$.}

\item \textit{If $\Sigma$ carries no nontrivial solution of  $\eqref{Bolle's eq}$ contractible in $M$,}
\beqn
\RFH(\Sigma,M)\cong\H(\Sigma;\Z/2).
\eeq
\end{enumerate}
In the extremal case, the assertions (i) and (iii) can be interpreted as:
\begin{itemize}
\item[(iv)] \textit{Let $T^n$ be a Lagrangian torus which is of restricted contact type with global coordinates embedded in $(M,\om)$ being convex at infinity. If an embedding $i:T^n\into M$ induces an injective homomorphism on $\pi_1$-level \footnote{ This implies that every solution of \eqref{Bolle's eq} is not contractible in $M$.} (e.g. the zero section of $T^*T^n$),
$$
\RFH(T^n,M)\cong\H(T^n;\Z/2).
$$
Therefore there always exists a self intersection point of $T^n$.}\\[-1.5ex]
\end{itemize}

\noindent\textbf{Further directions.} After finishing the present paper, various applications of Rabinowitz Floer theory have been conducted. We expect that most of those studies can be generalized to the coisotropic setting with the arguments in this paper.

\subsection{History and related results}\label{History and related results}

The aim of this paper is to extend Rabinowitz Floer homology theory to  coisotropic submanifolds. Rabinowitz Floer homology theory was developed by Cieliebak and Frauenfelder in \cite{CF} and has been extensively studied for hypersurfaces by many authors \cite{AF1,AF2,AF3,AF4,AF5,AM,CF,CFO,CFP,Ka1,Ka2,Ka3,AS,Me}. The framework and many results of this paper were inspired by their remarkable achievements and Ginzburg's pioneering work \cite{Gi1}.

The existence problem for leafwise intersection points was addressed by Moser \cite{M}. He obtained the result for simply connected symplectic manifolds and $C^1$-small perturbations and the simply connectedness assumption was removed by Banyaga \cite{Ba}. Hofer and Ekeland \cite{H,EH} replaced the assumption of $C^1$-smallness by the boundedness of the Hofer norm below a certain symplectic capacity for restricted contact hypersurfaces in $\mathbb{R}^{2n}$. Ginzburg \cite{Gi1} generalized Ekeland-Hofer's results for restricted contact hypersurfaces in subcritical Stein manifolds. Dragnev \cite{Dr} obtained this result on closed contact coisotropic submanifold in $\R^{2n}$. Albers-Frauenfelder \cite{AF1} proved the existence of leafwise intersection points for restricted contact hypersurfaces when the Hofer norm of perturbations is smaller than the minimal period of Reeb orbits. Using a different approach, G\"urel \cite{Gu} extended the theorem to restricted contact coisotropic submanifolds. Ziltener \cite{Z,Z2} also studied the question in a different way and obtained a lower bound of the number of leafwise intersection points under the assumption that the characteristic foliation is a fibration. The author \cite{Ka1} generalized Albers-Frauenfelder's theorem to unrestricted contact hypersurfaces but there was a constraint on the support of perturbations.

In a different aspect, the displacement energy of coisotropic submanifolds is also an integral part of the coisotropic intersection problem. Bolle \cite{Bo1,Bo2} proved that the displacement energy for stable coisotropic submanifolds of $\R^{2n}$ is positive. Ginzburg \cite{Gi1} extended Bolle's result to wide (or closed) and geometrically bounded manifolds. Recently, this was generalized further by Kerman \cite{Ke} and Usher \cite{U}.

\subsubsection*{Acknowledgments}
The author is deeply grateful to Urs Frauenfelder for many valuable discussions. He would like to thank Alex Oancea for sending the preprint \cite{BO2}. He also thanks the anonymous referees for numerous helpful comments on the manuscript. This work is partially supported by the Basic Research fund 2010-0007669.

\section{Coisotropic submanifolds}\label{coisotropic}
Let $(M,\om)$ be a $2n$-dimensional symplectic manifold. We denote by the {\em Hamiltonian tuple} $\GG:=(G_1,\dots,G_k)$ for time-independent Hamiltonian functions $G_i\in C^\infty(M)$, $i\in\{1,\dots,k\}$ for $1\leq k\leq n$. We often regard $\GG$ as an element of $C^\infty(M,\mathbb{R}^k)$.

\begin{Def}\label{def:Poisson-commute}
Given Hamiltonian functions $F$ and $G$ in $C^\infty(M)$, the {\em Poisson bracket}
$$
\{\cdot,\cdot\}:C^\infty(M)\x C^\infty(M)\pf C^\infty(M)
$$
is defined by $\{F,G\}:=\om(X_F,X_G)$.
A Hamiltonian tuple $\GG$ is said to be {\em Poisson-commute} if $\{G_i,G_j\}=0$ for all $1\leq i,j \leq k$.
\end{Def}

If a Hamiltonian tuple $\GG\in C^\infty(M,\R^k)$ Poisson-commutes and $0\in\R^k$ is a regular value of $\GG$, then $\GG^{-1}(0)$ is a smooth coisotropic submanifold of codimension $k$ in $(M,\om)$ and  $T\GG^{-1}(0)^\om$ is spanned by their Hamiltonian vector fields, namely $X_{G_1},\dots,X_{G_k}$. We refer to the introduction for the definitions of coisotropic, stable, contact, and restricted contact.\\[-1ex]

In the case that $\Sigma$ is of stable type, the normal bundle of $\Sigma\subset M$ is trivial, i.e. $N\Sigma\cong\Sigma\x\R^k$. From the Weinstein neighborhood theorem, we have
\begin{Prop}\cite{Bo1,Bo2}\label{prop:contact neighborhood}
Let $\Sigma$ be a closed stable coisotropic submanifold of codimension $k$ in $(M,\om)$. Then there exist $r>0$, a neighborhood $V$ of $\Sigma$ which is symplectomorphic by $\psi:U_r \to V$ to
$$
U_r:=\{(q,\pp)=(q,p_1,\dots,p_k)\in\Sigma\x\R^k\,|\,|p_i|<r,\,\,\textrm{for all }\, i=1,\dots,k\}
$$
with $\psi^*\om=\om_\Sigma+\sum_{i=1}^k d(p_i\alpha_i)$.
\end{Prop}

Here we use the same symbols $\om_\Sigma$ and $\alpha_i$ for differential forms in $\Sigma$ and for their pullback to $\Sigma\x\R^k$. We set
$$
\delta_0:=\max\bigr\{r\in\R\,\bigr|\,\textrm{there exists a symplectic embedding } \,\psi:U_r\into M\bigr\}
$$
and let $\psi_0:U_{\delta_0}\into M$ be a symplectic embedding. Henceforth, we identify $U_{\delta}$ with $\psi_0(U_{\delta})$ for all $0<\delta\leq\delta_0$. We have $X_{p_i}\in\ker\om_\Sigma$, $dp_j(X_{p_i})=0$ and $\alpha_j(X_{p_i})=\delta_{ij}$ on $\Sigma$ for $1\leq i,j\leq k$ since $i_{X_{p_i}}\om=dp_i$. Moreover the Hamiltonian tuple $\mathfrak{p}=(p_1,\dots,p_k)$ Poisson-commutes since $\{X_{p_1},\dots,X_{p_k}\}$ forms a basis for $\ker\om_\Sigma$. 

\begin{Def}
Let $\Sigma$ be a stable coisotropic submanifold in $(M,\om)$. 
The unique vector fields $R_1,\dots,R_k$ on $\Sigma$ characterized by  
$$
\alpha_i(R_j)=\delta_{ij},\quad R_i\in\bigcap_{\ell=1}^k\ker d\alpha_\ell,\quad i,\,j\in\{1,\dots,k\}
$$ 
are called the {\em Reeb vector fields} associated with the stable structure $(\Sigma,\alpha_1,\dots,\alpha_k)$.
\end{Def}

We note that $X_{p_1},\dots,X_{p_k}$ correspond to $R_1,\dots,R_k$ via the identification $\psi_0$. From now on, we choose an almost complex structure $J$ on $M$ which splits on $U_\epsilon$ with respect to
\beq\label{eq:spiltting of J}
TU_{\delta_0}=\bigg(\underbrace{\bigcap_{i=1}^k\ker\alpha_i}_{=:\xi}\bigg)\bigoplus\underbrace{\bigg(T\Sigma^\om\oplus\frac{\partial}
{\partial p_1}\oplus\cdots\oplus\frac{\partial}{\partial p_k}\bigg)}_{=:\xi^\om}.
\eeq
i.e. $J|_{\xi^\om}$ is an almost complex structure which interchanges the Reeb vector fields $R_i$ with $\frac{\partial}{\partial p_i}$ for $1\leq i\leq k$; strictly speaking $JR_i=\frac{\p}{\p p_i}$ and $J\frac{\p}{\p p_i}=-R_i$.

\begin{Def}\label{Def:global Bolle's coordinates}
We say that a closed coisotropic submanifold $\Sigma$ in $(M,\om)$ is of stable or (restricted) contact type {\em with global coordinates} if we are able to extend its coordinate functions $p_1,\dots,p_k$($\in C^\infty(U_{\delta_0})$ via $\psi_0$) to $\tilde p_1,\dots,\tilde p_k\in C^\infty(M)$ so that
\beqn
\tilde p_i:M\rightarrow\mathbb{R}\,\,\ by \,\,\,
  \tilde p_i = \left\{ \begin{array}{lll}
  p_i  & \textrm {on} &  U_{{\delta_0}/2}\\[1ex]
  g_i(p_i) & \textrm{on} & U_{\delta_0}\\[1ex]
  constant \qquad & \textrm  {outside} & U_{\delta_0}.\end{array}\right.\;
\eeq
for some $g_i\in C^\infty((-\delta_0,\delta_0))$ with $g_i^{-1}(0)=\{0\}$. We relabel $p_i$ instead of $\tilde p_i$ for notational convenience.
\end{Def}
\begin{Rmk}
It can easily be checked that the zero section of the cotangent bundle of a torus is of restricted contact type with global coordinates. Suppose that  $\Sigma_i$ in $(M_i,\om_i)$ for $i=1,2$ is a contact coisotropic submanifold with global coordinates. If in addition $\Sigma_1\x\Sigma_2$ is of contact type in $(M_1\x M_2,\om_1\oplus\om_2)$, it has global coordinates.

In order to define Rabinowitz Floer homology, we need global coordinates of stable or (restricted) contact coisotropic submanifolds. Therefore we encounter the problem under which conditions we are able to extend coordinates $p_1,\dots,p_k$ globally. This extending problem seems not easy, and that is why Step 2s in the proofs of Theorem A and Theorem D appear.
\end{Rmk}

\subsection{Examples of contact coisotropic submanifolds}
Although the contact condition is quite restrictive, we have the following  examples.
\begin{enumerate}
\item A coisotropic submanifold which is $C^1$-close to a contact coisotropic submanifold also has contact type.
\item A hypersurface has contact type if and only if it is of contact type in the standard sense.
\item A Lagrangian torus is of contact type with contact one forms $d\theta_1,\dots,d\theta_n$ where $\theta_1,\dots,\theta_n$ are angular coordinates on the $n$-dimensional torus. Indeed it turns out that a closed Lagrangian submanifold of contact type is necessarily a torus.
\item Let $\Sigma\subset (M_1,\om_1)$ be a contact coisotropic submanifold and $T^{n_2}\subset (M_2,\om_2)$ be a Lagrangian torus. Then a coisotropic submanifold $\Sigma\x T^{n_2}$ in $(M_1\x M_2,\om_1\oplus\om_2)$ is of contact type. In particular, the stabilization of $\Sigma\subset (M,\om)$, $\Sigma\x S^1\subset (M\x T^*S^1,\om\oplus d\theta\wedge dt)$ is of (restricted) contact type whenever $\Sigma$ is of (restricted) contact type. Here $\theta$ is the base coordinate and $t$ is the fiber coordinate.
\item Consider the Hopf fibration $\pi:S^{2n-1}\to\mathbb{C}P^{n-1}$. According to Marsden-Weinstein-Meyer reduction, we know that there is a canonical symplectic form $\om_{\mathbb{C}P^{n-1}}$ on $\mathbb{C}P^{n-1}$ satisfying $\pi^*\om_{\mathbb{C}P^{n-1}}=\om_{\R^{2n}}|_{S^{2n-1}}$ where $\om_{\R^{2n}}$ is the standard symplectic form on $\R^{2n}$. For a contact hypersurface $(\Delta,\alpha)\subset\mathbb{C}P^{n-1}$, $\pi^{-1}(\Delta)$ is a contact submanifold in $\R^{2n}$ of codimension $2$.
\end{enumerate}

Let $(M,\om)$ be a closed symplectic manifold with an integral symplectic form $[\om]\in\H^2(M;\Z)$. For each $N\in\N$, there exists a complex line bundle $p:E^{N}\to M$ with the first Chern class $c_1(E^{N})=-N[\om]$. We note that $S^1$ acts on the bundle $E^{N}$ by
\bean
S^1\x E^{N}&\pf E^{N}\\
(t,v)&\longmapsto e^{2\pi it}v.
\eea
Thus by the Boothby-Wang theorem, there exists a connection 1-form $\alpha$ on $E^{N}\setminus E_0$ where $E_0$ is the zero section of the complex line bundle $E^{N}\stackrel{p}{\to} M$; moreover it holds that  $p^*F_\alpha=d\alpha$ for the curvature 2-form $F_\alpha=N\om$. We  abbreviate $r=|e|$ for $e\in E^{N}$ and define $q:\R\to\R$ by $q(r)=\pi r^2+1/N$. Then the following two form gives a symplectic structure on $E^{N}$:
\beqn
\Omega_E:=q'(r)dr\wedge\alpha+q(r)N p^*\om.
\eeq
It is easy to check that $\Omega_E|_{E_0}=p_1^*\om$ and $\Omega_E|_{E\setminus E_0}=d(q(r)\alpha)$. Furthermore, for all $c>1/{N}$, the following submanifold
\beqn
\Sigma_c:=\{q(r)=c\}
\eeq
is of contact type. We perform this construction once again. We choose a complex line bundle $p':F^{K}{\to} M$ with the first Chern class $c_1(F^{K})=-K[\om]$. As before, there is a connection 1-form $\beta$ on $F^{K}\setminus F_0$ where $F_0$ is the zero section of the bundle $F^{K}\stackrel{p'}{\to} M$ such that its curvature 2-form $F_\beta$ satisfies $F_\beta=K\om$. We set the function $h(s)=\pi s^2+1/K$ for $s=|f|\in\R$ where $f\in F^K$, then
\beqn
\Omega_F:=h'(s)ds\wedge\beta+h(s){K} p'^*\om
\eeq
is a symplectic form on $F^{K}$. Next, we consider the Whitney sum of $E^N$ and $F^K$, $E^N\oplus F^K$ and let $\pi_1:E^N\oplus F^K\to E^N$ and $\pi_2:E^N\oplus F^K\to F^K$ be natural projections. We abbreviate $\tilde\om:=(p\circ\pi_1)^*\om=(p'\circ\pi_2)^*\om$, and use the same symbols $r$, $s$, $g(r)$, $h(s)$, $\alpha$, and $\beta$ for their pull-backs to $E^N\oplus F^K$. Then the following 2-form
\beqn
\Omega_{E\oplus F}:=h'(s)ds\wedge\beta+q'(r)dr\wedge\alpha +(q(r)N+h(s)K)\tilde \om
\eeq
becomes a symplectic form on $E^N\oplus F^K$. We have
\begin{itemize}
\item [(vi)] For any $c>1/{N}$ and $d>1/{K}$, set 
$$
\Delta_{c,d}:=\{q(r)=c,\,h(s)=d\}.
$$ 
Since $\Omega_{E\oplus F}|_{\Delta_{c,d}}=(cN+dK)\tilde\om$,  $\Delta_{c,d}$ with 1-forms $\frac{cN+dK}{N}\alpha$ and $\frac{cN+dK}{K}\beta$ is a contact coisotropic submanifold in $(E^N\oplus F^K,\Omega_{E\oplus F})$ of codimension 2.
\end{itemize}

\begin{Rmk}\cite{Bo2,Gi1}
Let $\Sigma$ be a closed contact coisotropic submanifold in $(M,\om)$. Then a 1-form $\lambda=a_1\lambda_1+\cdots+a_k\lambda_k$ with $a_1+\cdots+a_k=0$ is closed and represents an element in $\H^1_\mathrm{dR}(\Sigma)$. In addition, $\lambda\neq0$ is not exact; otherwise $\lambda=df$ for some $f\in C^1(\Sigma)$, $\lambda(x)=0$ at a critical point $x$ of $f$, but condition (ii) yields that $\lambda_1,\dots,\lambda_k$ are linearly independent on $\Sigma$; thus $\lambda_1(x)=\cdots\lambda_k(x)=0$. As a result, $\dim\H^1_\mathrm{dR}(\Sigma)\geq k-1$. It imposes restriction on the contact condition that a product of contact type coisotropic submanifolds is not necessarily of contact type; for instance, $S^3\x S^3$ is not of contact type in $\R^8$.
\end{Rmk}
\begin{Rmk}
Different from the contact case, a product of stable coisotropic submanifolds is of stable type again. Furthermore, a connected sum of a contact coisotropic submanifold is not of contact type in general; for instance, a connected sum of Lagrangian tori is not a torus any more, hence cannot be of contact type.
\end{Rmk}

\section{Rabinowitz action functional with several Lagrange multipliers}
From now on, our $2n$-dimensional symplectic manifold $(M,\om)$ is assumed to be symplectically aspherical, i.e. $\om|_{\pi_2(M)}=0$.

\subsection{Rabinowitz action functional}
Let $\eta=(\eta_1,\dots,\eta_k)\in\R^k$ be a $k$-tuple of Lagrange multipliers. Let $\LLL\subset C^\infty(S^1,M)$ be the component of contractible loops in $M$. For an arbitrary Poisson-commute Hamiltonian tuple $\GG=(G_1,\dots,G_k)\in C^\infty(M,\R^k)$ which has $0\in\R^k$ as a regular value, the generalized Rabinowitz action functional $\AA^\GG:\LLL\x\mathbb{R}^k\to\mathbb{R}$ is defined as follows:
\beq\label{eq:Rabinowitz action functional}
\AA^\GG(v,\eta):=-\int_{D^2}\bar v^*\om-\sum_{i=1}^k\eta_i\int_0^1G_i(v(t))dt
\eeq
where $\bar v$ is any filling disk of $v$, i.e. $\bar v|_{\partial D^2}(t)=v(t)$ for $t\in S^1$. The symplectic asphericity condition implies that the value of the above action functional is independent of the choice of filling discs. Using the standard scalar product $\langle\cdot,\cdot\rangle$ in $\mathbb{R}^k$, we can express \eqref{eq:Rabinowitz action functional} as
\beqn
\AA^\GG(v,\eta)=-\int_{D^2}\bar v^*\om-\int_0^1\langle\eta,\GG\rangle(v(t))dt.
\eeq

\subsection{Critical points}
A critical point of the Rabinowitz action functional, $(v,\eta)\in \Crit\AA^\GG$ satisfies the following equations.
\beq\label{eqn:critical point equation}\left.
\begin{aligned}
&\partial_tv(t)=\sum_{i=1}^k\eta_i X_{G_i}(v(t)),\quad  t\in S^1\\[1ex]
&\int_0^1G_i(v(t))dt=0, \quad   i\in\{1,\dots,k\}
\end{aligned}
\;\;\right\}
\eeq

\begin{Prop}\label{Prop:critical points lie in contact}
If $(v,\eta)\in \Crit\AA^\GG$, $v$ lies in the coisotropic submanifold $\GG^{-1}(0)$.
\end{Prop}
\begin{proof}
Using the first equation in \eqref{eqn:critical point equation}, we have
\beqn
\frac{d}{dt} G_i(v(t))=dG_i(v(t))[\partial_t v]=dG_i\bigg(\sum_{j=1}^k\eta_j X_{G_j}(v(t))\bigg)=\sum_{j=1}^k\eta_j\underbrace{\{G_i,G_j\}}_{=0}(v(t))=0
\eeq
It means that every Hamiltonian function $G_i$ is stationary along $v(t)$. The second equation in \eqref{eqn:critical point equation} implies $G_i(v(t))=0$ for all $1\leq i\leq k$, $t\in S^1$. This proves the proposition.
\end{proof}

Let $\chi\in C^{\infty}(S^1,\mathbb{R})$ be a smooth function such that $\chi(t)\geq 0$, $\int_0^1\chi(t)dt=1$, and $\Supp\chi\subset(1/2,1)$. Suppose that $\Sigma$ is a stable coisotropic submanifold with  global coordinates  $\pp=(p_1,\dots,p_k)\in C^\infty(M,\R^k)$. Using $\chi$ with $\pp$, we define a time-dependent Hamiltonian $H_i:S^1\x M\rightarrow\mathbb{R}$ by $H_i(t,x)=\chi(t)p_i(x)$ for all $1\leq i\leq k$, i.e.
\beq\label{eq:defining Hamiltonian}
\HH(t,x):=\chi(t)\pp(x)\in C^\infty(S^1\x M,\R^k).
\eeq
We remark that the Rabinowitz action functional can be defined with an arbitrary Hamiltonian tuple and Proposition \ref{Prop:critical points lie in contact} holds for an arbitrary Poisson-commuting Hamiltonian tuple; but from now on, we only deal with the Hamiltonian tuple $\HH$.

In the restricted contact hypersurface case \cite{CF}, we know that if $(v,\eta)\in\LLL\x\R$ is a critical point of the Rabinowitz action functional then $v_\eta(t):=v(t/\eta)$ is a closed Reeb orbit lying on a hypersurface with period $\eta$. In the similar vein, we have

\begin{Lemma}\label{Lem:eta's are same}
Let $\Sigma$ be a restricted contact coisotropic submanifold with global coordinates in $(M,\om)$. If $(v,\eta)=(v,\eta_1,\dots,\eta_k)\in \Crit\AA^\HH$, then $\eta_1=\cdots=\eta_k$.
\end{Lemma}
Thus if $(v,\eta)\in\Crit\AA^\HH$ with $\eta\neq0$, then $v_\eta(t):=v(t/\eta_1)$ is a solution of \eqref{equation}, i.e.
\beq\label{eq:eta is period}
\partial_t v_\eta(t)=\sum_{i=1}^k  R_i(v_\eta(t)),\quad t\in \R/\eta\Z.
\eeq
When $\eta=0$, $v$ is a constant loop in $\Sigma$, i.e. $\p_tv\equiv 0$.

\begin{proof}(of Lemma \ref{Lem:eta's are same})
Due to the previous proposition, $v(t)\in\Sigma$ for all $t\in S^1$. Since all $\lambda_i$'s are primitives of $\om$, we have
\bean
\AA^\HH(v,\eta)
&=-\int_{D^2}\bar v^*\om-\int_0^1\langle\eta,\HH\rangle(t,v(t)) dt\\
&=-\int_0^1v^*\lambda_i-\sum_{j=1}^k\eta_j\int_0^1\chi(t) \underbrace{p_i(v(t))}_{0}dt\\
&=-\int_0^1\alpha_i\Big(\sum_{j=1}^k\eta_jX_{H_j}(v)\Big)dt\\
&=-\eta_i\int_0^1\alpha_i(X_{H_i}(v))dt\\
&=-\eta_i
\eea
for all $i\in\{1,\dots,k\}$ where $\lambda_i|_\Sigma=\alpha_i$. Thus we conclude that all $\eta_i$'s coincide.
\end{proof}

\section{Perturbation of the generalized Rabinowitz action functional}
\subsection{Hofer norm}
Let $\Ham_c(M,\om)$ be the group of Hamiltonian diffeomorphisms generated by compactly supported Hamiltonian functions. We briefly recall the definition of {\em Hofer norm}.
\begin{Def}\label{Def:Hofer norm}
Let $F\in C_c^\infty(S^1\x M,\mathbb{R})$ be a compactly supported Hamiltonian function. We set
$$
||F||_+:=\int_0^1\max_{x\in M} F(t,x) dt,\qquad||F||_-:=-\int_0^1\min_{x\in M} F(t,x) dt=||-F||_+.
$$
The Hofer norm of $F$ is defined by
$$
||F||=||F||_++||F||_-.\;
$$
For $\phi\in \Ham_c(M,\om)$, the Hofer norm is
$$
||\phi||=\inf\{||F||\mid \phi=\phi_F,\,F\in C^\infty_c(S^1\x M,\R)\}.\;
$$

\end{Def}

\begin{Lemma}\cite{AF1}\label{lemma:norms equivalent}
For all $\phi\in\Ham_c(M,\om)$,
$$
||\phi||=|||\phi|||:=\inf\{||F||\mid \phi=\phi_F,\;F(t,\cdot)=0\;\;\forall t\in[\tfrac12,1]\}\;.
$$
\end{Lemma}
\begin{proof}
To prove $||\phi||\geq|||\phi|||$, pick a smooth monotone increasing map $r:[0,1]\to[0,1]$ with $r(0)=0$ and $r(\tfrac12)=1$. For $F$ with $\phi_F=\phi$ we set $F^r(t,x):=r'(t)F(r(t),x)$. Then a direct computation shows $\phi_{F^r}=\phi_F$, $||F^r||=||F||$, and $F^r(t,x)=0$ for all $t\in[\tfrac12,1]$. The reverse inequality is obvious.
\end{proof}

\subsection{Perturbed Rabinowitz action functional}
In this subsection $\Sigma$ is a closed stable coisotropic submanifold in $(M,\om)$. Let $\HH\in C^\infty(S^1\x M,\R^k)$ be a Hamiltonian tuple as in \eqref{eq:defining Hamiltonian} and $F\in C_c^\infty(S^1\x M)$ be an arbitrary time-dependent Hamiltonian function. Thanks to the previous lemma we assume that $F$ has time support in $(0,\frac{1}{2})$. We note that the time support of $\HH$ and the time support of $F$ are {\em disjoint}. With these Hamiltonian functions, the perturbed Rabinowitz action functional $\AA^\HH_F:\LLL\x\R^k\to\R$ is defined by
\beqn
\AA^\HH_F(v,\eta):=-\int_{D^2}\bar v^*\om-\int_0^1F(t,v(t))dt-\int_0^1\langle\eta,\HH\rangle(t,v(t))dt.
\eeq
where $\bar v:D^2\to M$ is any filling disk of $v$. A critical point of the perturbed Rabinowitz action functional, $(v,\eta)\in \Crit\AA^\HH_F$ satisfies the following equations.
\beq\label{eqn:critical point equation 2}\left.
\begin{aligned}
&\partial_tv(t)=X_F(t,v)+\sum_{i=1}^k\eta_i X_{H_i}(t,v(t)),\quad  t\in S^1\\[1ex]
&\int_0^1H_i(t,v(t))dt=0, \quad  i\in\{1,\dots,k\}
\end{aligned}
\;\;\right\}
\eeq

In the next proposition, we observe that a critical point of $\AA^\HH_F$ gives rise to a leafwise intersection point. Albers-Frauenfelder \cite{AF1} proved the following proposition when $\Sigma$ is a hypersurface. Following through their proof, we prove

\begin{Prop}\label{prop:critical point answers question}
If $(v,\eta)\in\Crit\AA_F^\HH$, $v(0)\in\Sigma$ is a leafwise intersection point.
\end{Prop}
\begin{proof}
Since the time support of $F$ is $(0,1/2)$, for $t\geq1/2$ and for all $i=1,\dots,k$,
\beqn
\frac{d}{dt}p_i(v(t))=dp_i(v(t))[\partial_t v]=dp_i(v(t))\bigr[\underbrace{X_F(t,v)}_{=0}+\sum_{j=1}^k\eta_j X_{H_j}(t,v)\bigr]=0
\eeq
The last equality follows from Poisson-commutativity of $\pp$. As in the proof of Proposition \ref{Prop:critical points lie in contact}, the second equation in \eqref{eqn:critical point equation 2} implies $v(t)\in \pp^{-1}(0)=\Sigma$ for $t\in(1/2,1)$. On the other hand, $v$ solves $\partial_t v=X_F(t,v)$ on $(0,1/2)$ so that $v(1/2)=\phi_F^{1/2}(v(0))=\phi_F^1(v(0))$ since $F=0$ for $t\geq1/2$. For $t\in(1/2,1)$, it holds that $\partial_t v=\sum_{i=1}^k\eta_i X_{H_i}(t,v)$ and thus $v(0)=v(1)\in L_{v(1/2)}$. Thus we conclude that $v(0)\in L_{\phi_F(v(0))}$ which is equivalent to $\phi_F(v(0))\in L_{v(0)}$.
\end{proof}

From now on, we allow $s$-dependence on $F$ as follows. Let $\{F_s\}_{s\in\mathbb{R}}$ be a family of Hamiltonian functions varying only on a finite interval in $\mathbb{R}$. More specifically, we assume $F_s(t,x) =F_-(t,x)$ for $s\leq-1$ and $F_s(t,x)=F_+(t,x)$ for $s\geq1$. We also choose a family of compatible almost complex structures $\{J(s,t)\}_{(s,t)\in{\R\x S^1}}$ on $M$ such that $J(s,t)$ is invariant outside of the interval $[-1,1]\subset\R$ and they still split as in \eqref{eq:spiltting of J}.

On the tangent space $T_{(v,\eta)}(\LLL\times\mathbb{R}^k)=T_v\LLL\x T_\eta\mathbb{R}^k$ for $(v,\eta)\in\LLL\x\R^k$, we define the metric $m$ as follows:
\beqn
m_{(v,\eta)}\big((\hat v^1,\hat\eta^1),(\hat v^2,\hat\eta^2)\big):=\int_0^1 g_{v}(\hat v^1,\hat v^2)dt+\langle\hat\eta^1,\hat\eta^2\rangle.
\eeq
where $g(\cdot,\cdot):=\om(\cdot,J\cdot)$ is the metric on $M$. Here $\hat\eta^1$ and $\hat\eta^2$ are elements in $T_\eta\R^k\cong\R^k$ and $\langle\cdot,\cdot\rangle$ is the scalar product in $\R^k$.
\begin{Def}
A map $w\in C^\infty(\mathbb{R},\LLL\times\mathbb{R}^k)$ which solves
\beq\label{eqn:gradient flow line}
\partial_s w(s)+\nabla_m \AA_{F_s}^\HH(w(s))=0.\;
\eeq
is called a {\em gradient flow line} of $\AA_{F_s}^\HH$ with respect to the metric $m$.
\end{Def}

According to Floer's interpretation, the gradient flow equation \eqref{eqn:gradient flow line} can be interpreted as $w=(u,\tau)=(u,\tau_1,\dots,\tau_k)$ with $u(s,t):\mathbb{R}\x S^1\to M$ and $\tau_i(s):\mathbb{R}\to\mathbb{R}$,  solving
\beq\label{eq:Floer's interpretation}\left.
\begin{aligned}
&\partial_su+J(s,t,u)\bigr(\partial_tu-\sum_{i=1}^k \tau_i X_{H_i}(t,u)-X_{F_s}(t,u)\bigr)=0\\[1ex]
&\partial_s\tau_i-\int_0^1H_i(t,u)dt=0, \qquad 1\leq i\leq k.
\end{aligned}
\;\;\right\}
\eeq

\begin{Def}
The energy of a map $w\in C^{\infty}(\R,\LLL\times\R^k)$ is defined as
\beqn
E(w):=\int_{-\infty}^\infty||\partial_s w||_m^2ds\;.
\eeq
\end{Def}

\begin{Lemma}\label{lemma:energy estimate for gradient lines}
Let $w\in C^\infty(\R,\LLL\x\R^k)$ be a gradient flow line of $\AA_{F_s}^\HH$ with finite energy. Then we have the following estimation.
\beq\label{eqn:energy estimate for gradient lines}
E(w)\leq\AA_{F_-}^\HH(w_-)-\AA_{F_+}^\HH(w_+)+\int_{-\infty}^\infty||\partial_sF_s||_-ds\;
\eeq
where $w_\pm:=\lim_{s\to\pm\infty}w(s)\in\Crit\AA^\HH_{F_s}$.
Moreover, equality holds if $\partial_sF_s=0$.                                                                                                                                                                                                                                                                                                                                                                                                                                                                                                                                                                                                                                                                                                                                                                                                                                                                                                                                                                                                                                                                                                                                                                                                                                                                                                                                                                                                                                                                                                                                                                                                                                                                                                                                                                                                                                                                                                                                                                                                                                                                                                                                                                                                                                                                                                                                                                                                                                                                                                                                                                                                                                                                                                                                                                                                                                                                                                                                                                                                                                                                                                                                                                                                                                                                                                                                                                                                                                                                                                                                                                                                                                                                                                                                                                                                                                                                                                                                                                                                                                                                                                                                                                                                                                                                                                                                                                                                                                                                                                                                                                                                                                                                                                                                                                                                                                                                                                                                                                                                                                                                                                                                                                                                                                                                                                                                                                                                                                                                                                                                                                                                                                                                                                                                                                                                                                                                                                                                                                                                                                                                                                                                                                                                                                                                                                                                                                                                                                                                                                                                                                                                                                                                                                                                                                                                                                                                                                                                                                                                                                                                                                                                                                                                                                                                                                                                                                                                                                                                                                                                                                                                                                                                                                                                                                                                                                                                                                                                                                                                                                                                                                                                                                                                                                                                                                                                                                                                                     \end{Lemma}

\begin{proof}
Let us consider $w=(u,\tau)\in C^\infty(\R\x S^1,M)\x C^\infty(\R,\R^k)$ as \eqref{eq:Floer's interpretation}.
\bean
E(u,\tau)&=-\int_{-\infty}^\infty d\AA_{F_s}^\HH\big((u,\tau)(s)\big)[\partial_s(u,\tau)]ds\\[1ex]
&=-\int_{-\infty}^\infty \frac{d}{ds}\Big(\AA_{F_s}^\HH\big((u,\tau)(s)\big)\Big)ds+\int_{-\infty}^\infty \big(\partial_s\AA_{F_s}^\HH\big)\big((u,\tau)(s)\big)ds\\[1ex]
&=\AA_{F_-}^\HH(w_-)-\AA_{F_+}^\HH(w_+)-\int_{-\infty}^\infty\int_0^1\partial_sF_s(t,u)dtds\\[1ex]
&\leq\AA_{F_-}^\HH(w_-)-\AA_{F_+}^\HH(w_+)+\int_{-\infty}^\infty||\partial_sF_s||_-ds\;.
\eea
This computation proves the lemma.
\end{proof}

\begin{Rmk}
We note that $\int_{-\infty}^\infty||\p_s F_s||_-ds$ has a finite value since $\p_s F_s$ has a compact support by construction.
\end{Rmk}

\begin{Prop}\label{prop:uniform bound of action functional}
$\AA_{F_s}^\HH$ has a uniform bound along gradient flow lines.
\end{Prop}
\begin{proof}
For any gradient flow line $w\in C^\infty(\R,\LLL\x\R^k)$ of $\AA_F^\HH$ and $s_1<s_2\in\R$, we calculate
\bean
 0&\leq\int_{s_1}^{s_2}||\p_s w||_m^2\,ds\\
 &=-\int_{s_1}^{s_2}d\AA_{F_s}^\HH(w(s))(\partial_s w)ds\\
 &=\AA_{F_{s_1}}^{\HH}(w(s_1))-\AA_{F_{s_2}}^{\HH}(w(s_2)) - \int_{s_1}^{s_2}\int_0^1\partial_s F_s(t,v)dtds\\
&\leq\AA_{F_{s_1}}^{\HH}(w(s_1))-\AA_{F_{s_2}}^{\HH}(w(s_2))+\int_{s_1}^{s_2}||\partial_s F_s||_-ds.\\
\eea
From the above inequality we obtain
\bean
&\AA_{F_{s_2}}^{\HH}(w(s_2))\leq\AA_{F_-}^{\HH}(w_-)+\int_{-\infty}^{\infty}||\partial_s F_s||_-ds,\\
&\AA_{F_{s_1}}^{\HH}(w(s_1))\geq\AA_{F_+}^{\HH}(w_+)-\int_{-\infty}^{\infty}||\partial_s F_s||_-ds.
\eea
This proves the proposition.
\end{proof}

\subsection{Compactness}
In this subsection, we prove Theorem \ref{thm:compactness of moduli} which is a vital ingredient for all our results. Here, $\Sigma$ is assumed to be a closed restricted contact coisotropic submanifold and $(M,\om)$ is convex at infinity. But if we impose some restriction on the perturbation $F$, the theorem holds true in the contact case as well, see the next section.
\begin{Thm}\label{thm:compactness of moduli}
Let $\{w^\nu=(u^\nu,\tau^\nu)\}_{\nu\in\N}$ be a sequence of gradient flow lines of  $\AA_{F_s}^\HH$ with a uniform action bound like \eqref{eq:assumption}.
Then for every reparametrization sequence $\sigma_\nu\in\R$ the sequence $w^\nu(\cdot+\sigma_\nu)$ has a convergent subsequence in the $C^\infty_{loc}$-topology. That is, $\{w^\nu\}_{\nu\in\N}$ has a subsequence which converges with all derivatives on every compact subset to a gradient flow line $w\in C^\infty(\R\x S^1,M)\x C^\infty(\R,\R^k)$.
\end{Thm}
\begin{proof}
Once we establish the following three ingredients,
\begin{enumerate}
\item a uniform $L^\infty$-bound on $u^\nu$;
\item a uniform $L^\infty$-bound on $\tau^\nu$;
\item a uniform $L^\infty$-bound on the derivatives of $u^\nu$,
\end{enumerate}
the proof of the theorem follows from the elliptic bootstrapping argument in Floer theory, see \cite[Theorem B.4.2]{MS}. (i) follows since the image of every $u^\nu$, $\nu\in\N$ lies in a compact subset of $M$ due to the assumption ``convex at infinity", see \cite[Lemma 2.4]{Mc}. (ii) is the new feature of Rabinowitz Floer theory and is proved in Theorem \ref{thm:bound on eta}. Then if the derivatives of $u^\nu$ explode, the ``bubbling-off'' phenomenon occurs, see \cite[Chapter 4.2]{MS}. That is, we can detect non-constant $J$-holomorphic spheres as limits; but the symplectic asphericity of $(M,\om)$ rules out this possibility. Hence (iii) and the theorem is proved.
\end{proof}

We recall that $\lambda_i(X_{p_j})|_\Sigma=\delta_{ij}$, see Proposition \ref{prop:contact neighborhood}. Thus we can pick $\delta_1\in(0,\delta_0)$, so that on $U_{\delta_1}=\mathfrak{p}^{-1}(-\delta_1,\delta_1)$ and for all $1\leq i\neq j\leq k$,
\bea\label{eq:near contact}
&\frac{3}{4}<\lambda_i(X_{p_i})<\frac{5}{4}, \qquad &-\frac{1}{(4k-1)}<\lambda_i(X_{p_j})<\frac{1}{(4k-1)}.
\eea
Then we can prove the following fundamental lemma.
\begin{Lemma}\label{lemma:bound on eta when gradient is small}
There exist $\epsilon>0$ and $C>0$ such that for $(v,\eta)\in\LLL\x\mathbb{R}^k$,
\beqn
||\nabla_m\AA^\HH_{F_s}(v,\eta)||_m<\epsilon \;\;\textrm{  implies  }\;\;|\eta_i|\leq C\big(|\AA^\HH_{F_s}(v,\eta)|+1\big) \;\textrm { for all } 1\leq i\leq k.
\eeq
\end{Lemma}

\begin{proof}
The proof proceeds in three steps.\\[0.5ex]
\textbf{Step 1:} Assume $v(t)\in U_\delta$ for $t\in (1/2,1)$ and $\delta=\min\{1/4k,\delta_0/2,\delta_1\}$. Then there exists $C_0>0$ satisfying the following inequality for all $i=1,\dots,k$:
\beqn
|\eta_i|\leq C_0\big(|\AA^\HH_{F_s}(v,\eta)|+||\nabla_m\AA^\HH_{F_s}(v,\eta)||_m+1\big).
\eeq
Proof of Step 1. For each $i=1,\dots,k$, we estimate,
\bean
|\AA^\HH_{F_s}(v,\eta)|
=&\bigg|\int_0^1v^*\lambda_i+\sum_{j=1}^k\eta_j\int_0^1H_j(t,v)dt+\int_0^1F_s(t,v)dt\bigg|\\
\geq&\bigg|\sum_{j=1}^k\eta_j\int_0^1\lambda_i(v)\big(X_{H_j}(t,v)\big)dt\bigg|-\bigg|\int_0^1\lambda_i(v)\big(X_{F_s}(t,v)\big)dt\bigg|-\bigg|\int_0^1F_s(t,v)dt\bigg|\\
&-\bigg|\int_0^1\lambda_i(v)\big(\partial_t v-\sum_{j=1}^k\eta_jX_{H_j}(t,v)-X_{F_s}(t,v)\big)dt\bigg|-\bigg|\sum_{j=1}^k\eta_j\int_0^1H_j(t,v)dt\bigg|\\
\geq&\bigg|\eta_i\int_{1/2}^1\lambda_i(v)\big(X_{H_i}(t,v)\big)\bigg|-\bigg|\sum_{j\ne i}\eta_j\int_{1/2}^1\lambda_i(v)\big(X_{H_j}(t,v)\big)\bigg|-\bigg|\int_0^1\lambda_i(v)\big(X_{F_s}(t,v)\big)dt\bigg|\\
&-C_i||\nabla_m\AA^\HH_{F_s}(v,\eta)||_m-\bigg|\sum_{j=1}^k\eta_j\int_{1/2}^1H_j(t,v)dt\bigg|-\bigg|\int_0^1F_s(t,v)dt\bigg|\\
\geq&\frac{3}{4}|\eta_i|-\frac{1}{4(k-1)}\sum_{j\ne i}|\eta_j|-C_i||\nabla_m\AA^\HH_{F_s}(v,\eta)||_m-\delta\sum_{j=1}^k|\eta_j|-C_{i,F}
\eea
where $C_i:=||\lambda_i|_{U_\delta}||_{L^\infty}<\infty$ and $C_{i,F}:=||F_s||_{L^\infty}+C_i||X_{F_s}||_{L^\infty}<\infty$ for $i=1,\dots,k$. From the above inequality for all $i=1,\dots,k$, we extract the following estimation.
\beqn
\frac{1}{4}\sum_{i=1}^k|\eta_i|\leq k|\AA^\HH_{F_s}(v,\eta)|+\sum_{i=1}^k\big(C_i||\nabla_m\AA^\HH_{F_s}(v,\eta)||_m+C_{i,F}\big).
\eeq
Therefore Step 1 follows with $C_0:=\max\{4k,4kC_1,\dots,4kC_k,4kC_{1,F},\dots,4kC_{k,F}\}$.\\

\noindent\textbf{Step 2:} If there exists $t\in (\frac{1}{2},1)$ such that $v(t)\notin U_\delta$ then $||\nabla_m \AA^\HH_{F_s}(v,\eta)||_m\geq\epsilon$.\\[1ex]
\noindent Proof of Step 2. The assumption $v(t)\notin U_\delta$ means that there exists $i\in\{1,\dots,k\}$ such that $v(t)\notin U_\delta^i=p_i^{-1}(-\delta,\delta)$.
If in addition, $v(t) \in M-U^i_{\delta/2}$ for all $t \in (\frac{1}{2},1)$ then we easily have
\beqn
||\nabla_m\AA_{F_s}^\HH(v,\eta)||_m\geq\bigg|\int_0^1H_i(t,v(t))dt\bigg|=\bigg|\int_{1/2}^1\chi(t)p_i(v(t))dt\bigg|\geq\frac{\delta}{2}.\;
\eeq
Otherwise there exists $t' \in(\frac{1}{2},1)$ such that $v(t')\in U^i_{\delta/2}$. Thus we can find $t_0, t_1 \in (\frac{1}{2},1)$ such that
\bean
&v(t_0)\in \partial U^i_{\delta/2}, v(t_1)\in \partial U^i_{\delta}
\quad\&\quad v(t)\in U^i_\delta-U^i_{\delta/2},\quad \textrm{for  }\,\forall t\in [t_0,t_1],\\
\textrm{or}\qquad&
 v(t_1)\in \partial U^i_{\delta}, v(t_0)\in \partial U^i_{\delta/2}
\quad\&\quad v(t)\in U^i_\delta-U^i_{\delta/2},\quad \textrm{for  }\,\forall t\in [t_1,t_0].
\eea
We treat only the first case. The latter case is analogous. With $\mathfrak{P}:=\max_{x\in U_{\delta}}||\nabla_g p_i(x)||_g<\infty$ we estimate,
\bea\label{eq:main ineq in step 2 in fundamental lemma}
\mathfrak{P}||\nabla_m\AA^\HH_{F_s}(v,\eta)||_m
&\geq \mathfrak{P}||\partial_t v-\sum_{j=1}^k\eta_j X_{H_j}(t,v)-X_{F_s}(t,v)||_{L^2}\\
&\geq \mathfrak{P}||\partial_t v-\sum_{j=1}^k\eta_j X_{H_j}(t,v)-X_{F_s}(t,v)||_{L^1}\\
&\geq \int_{t_0}^{t_1}||\partial_t v-\sum_{j=1}^k\eta_j X_{H_j}(t,v)-X_{F_s}(t,v)||_g||\nabla_g p_i(v(t))||_g dt\\
&\geq \bigg|\int_{t_0}^{t_1}\big\langle\nabla_g p_i(v(t)),\partial_t v(t)-\sum_{j=1}^k\eta_j X_{H_j}(t,v)-X_{F_s}(t,v)\big\rangle_g dt\bigg|\\
&= \bigg|\int_{t_0}^{t_1}dp_i(v(t))\bigr(\partial_t v(t)-\sum_{j=1}^k\eta_j X_{H_j}(t,v)-\underbrace{X_{F_s}(t,v)}_{=0})\bigr) dt\bigg|\\
&= \bigg|\int_{t_0}^{t_1}\frac{d}{dt} p_i(v(t))dt-\underbrace{dp_i(v)\big(\sum_{j=1}^k\eta_j X_{H_j}(t,v)\big)}_{=0}\bigg|\\
&\geq \big|p_i(v(t_1))\big|-\big|p_i(v(t_0))\big|\\
&=\frac{\delta}{2}.
\eea
Thus Step 2 follows with $\epsilon=\min\bigr\{\frac{\delta}{2},\frac{\delta}{2\mathfrak{P}}\bigr\}.$\\[-1ex]

\noindent\textbf{Step\, 3:} Proof of the lemma.\\[1ex]
Proof of Step 3. According to Step 2, $v(t)\in U_\delta$ for all $t\in (\frac{1}{2},1).$ Then Step 1 completes the proof of the lemma with $C=C_0+\epsilon+1.$
\end{proof}
For a given gradient flow line $w$ of $\AA^\HH_{F_s}$ and $\sigma\in\mathbb{R}$, we define
\bea\label{eq:def of tau and C_F}
o(\sigma,w,\epsilon)&:=\inf\bigr\{\tau\geq0\,\bigr|\,\,||\nabla_m\AA^\HH_{F_s}(w(\sigma+\tau))||_m\leq\epsilon\bigr\},\\
C_F&:=\int_{-\infty}^{\infty}\int_0^1\max_{x\in M}||\partial_s F_s(t,x)||_g dtds \,\,<\,\, \infty.
\eea
\begin{Lemma}\label{lemma:tau(sigma)}
For a gradient flow line $w$ of $\AA^\HH_{F_s}$ with $\lim_{s\to\pm\infty}w(s)=w_\pm$,
\beqn
o(\sigma,w,\epsilon)\leq\frac{\AA^\HH_{F_s}(w_-)-\AA^\HH_{F_s}(w_+)+C_F}{\epsilon^2}.
\eeq
\end{Lemma}
\begin{proof}
We compute
\bean
\epsilon^2o(\sigma,w,\epsilon)
&\leq\int_\sigma^{\sigma+o(\sigma,w,\epsilon)}\big|\big|\nabla_m \AA_{F_s}^\HH(w)\big|\big|_m^2ds\\
&\leq\int_{-\infty}^{\infty}-d\AA_{F_s}^\HH(w)(\partial_s w)ds-C_F+C_F\\
&\leq\int_{-\infty}^{\infty}-\frac{d}{ds}\bigr(\AA_{F_s}^\HH(w(s))\bigr)ds+C_F\\
&=\AA_{F_s}^\HH(w_-)-\AA_{F_s}^\HH(w_+)+C_F
\eea
We obtain a bound on $o(\sigma,w,\epsilon)$ by dividing $\epsilon^2$ in the above inequality.
\end{proof}

\begin{Thm}\label{thm:bound on eta}
Assume that $w=(u,\tau)\in C^{\infty}(\mathbb{R},\LLL\x\mathbb{R}^k)$ is a gradient flow line of $\AA_{F_s}^\HH$ for which there exist $a \leq b$ such that
\beq\label{eq:assumption}
a\leq\AA_{F_s}^\HH(w(s))\leq b, \quad\textrm{for all } s\in \mathbb{R}.
\eeq
Then the $L^\infty$-norms of $\tau_i$'s are uniformly bounded.
\end{Thm}
As we have mentioned, Theorem \ref{thm:bound on eta} completes the proof of Theorem \ref{thm:compactness of moduli}.
\begin{proof}
Using Lemma \ref{lemma:bound on eta when gradient is small} and Lemma \ref{lemma:tau(sigma)}, we obtain
\bean
|\tau_i(\sigma)|&\leq|\tau_i(\sigma+o(\sigma,w,\epsilon))|+\int_\sigma^{\sigma+o(\sigma,w,\epsilon)}|\partial_s\tau_i(s)|ds\\
&\leq C\bigr(\bigr|{\AA}^\HH_{F_s}(w(\sigma+o(\sigma,w,\epsilon)))\bigr|+1\bigr)+o(\sigma,w,\epsilon)||H_i||_{L^\infty}\\
&\leq C(\max\{|a|,|b|\}+1)+
\Bigg(\frac{|b-a|+C_F}{\epsilon^2}\Bigg)||H_i||_{L^\infty}.
\eea
\end{proof}

\subsection{Proof of Theorem A}
The proof proceeds in two steps. In Step 1, we prove Theorem A under the assumption that $\Sigma$ admits global coordinates (separating condition). Then we remove this additional assumption in Step 2.\\[0.5ex]
\noindent\textbf{Step 1.} There exists a critical point $(v,\eta)$ of $\AA_F^\HH$ if $||F||<\wp(\Sigma)$ and $\Sigma$ is of restricted contact type with global coordinates. Moreover the action value of that critical point is uniformly bounded as below:
\beq\label{eq:step 1 in thm A}
-||F||\leq\AA_F^\HH(v,\eta)\leq||F||.
\eeq

\noindent{Proof of Step 1}. We mainly follow the proof of Theorem A in \cite{AF1} which made use of the ``stretching the neck'' argument. For $0\leq r$, we choose a smooth family of functions $\varphi_r\in C^{\infty}(\mathbb{R},[0,1])$ satisfying
\begin{enumerate}
\item for $r\geq1$: $\varphi_r'(s)\cdot s\leq0$ for all $s\in\R$, $\varphi_r(s)=1$ for $|s|\leq r-1$, and $\varphi_r(s)=0$ for $|s|\geq r$,
\item for $r\leq1$: $\varphi_r(s)\leq r$ for all $s\in\R$ and $\Supp\varphi_r\subset[-1,1]$,
\end{enumerate}

We note that $\varphi_\infty\equiv 1$ is the limit of $\varphi_r$ with respect to $C^\infty_{loc}$-topology. We fix a point $p\in\Sigma$ and consider the moduli space
$$
\M:=\left\{(r,w)\in[0,\infty)\times C^{\infty}(\R,\LLL\times\R^k)\,\bigg|\, \begin{aligned}& \,\,\,w\text{ is a gradient flow line of $\AA_{\varphi_rF}^\HH$ with}\\&\!\lim_{s\to-\infty}w(s)=(p,0),\;\lim_{s\to\infty}w(s)\in\Sigma\x\{0\} \end{aligned}\right\}\;.
$$
Assume on the contrary that there is {\em no} leafwise intersection point of $\phi_F$ for $||F||<\wp(\Sigma)$. For $(r,w)\in\M$ with $w_-=(p,0)$ and $w_+=(q,0)$ in $\Sigma\x\{0\}$, we estimate
\bean
E(w)&=-\int_{-\infty}^{\infty}d\AA_{\varphi_r(s)F}^\HH(w(s))(\partial_s w) ds\\
&\leq\AA_0^\HH(p,0)-\AA_0^\HH(q,0)+\int_{-\infty}^\infty||\partial_s\varphi_r F||_-ds\\
&=\int_{-\infty}^\infty||\varphi_r'(s)F||_-ds\\
&=\int_{-\infty}^0\varphi_r'(s)||F||_-ds-\int_{0}^\infty\varphi_r'(s)||F||_+ds\\
&=\varphi_r(0)\big(||F||_-+||F||_+\big)\\
&\leq||F||.
\eea
Accordingly we can also estimate,
\beq\label{eq:uniform action bounds}
-||F||\leq\AA_{\varphi_{r_n}F}^\HH(w_n(s))\leq||F||,\qquad (r_n,w_n)\in\M.
\eeq
Due to the action bound, Theorem \ref{thm:compactness of moduli} yields that a sequence $\{w_n\}_{n\in\N}$ for $(r_n,w_n)\in\M$ has a convergent subsequence (still denoted $w_n$) in $C^\infty_{loc}$-topology. We denote by $x$ the limit gradient flow line (which can be a constant gradient flow line). We want to show that $\M$ is compact and so assume by contradiction that $x_+\notin\Sigma\x\{0\}$ where $x_\pm$ are asymptotic ends of $x$, i.e. $x_\pm=\lim_{s\to\pm\infty}x(s)$.\\[-1.5ex]

\noindent\underline{Case 1}. $r_n$ is bounded.\\[-1.5ex]

There is no loss of generality in assuming that $r_n\to r$ as $n\to\infty$. Let $U\in\LLL\x\R^k$ be an open set containing only the constant critical points of $\AA^\HH_{\varphi_r F}$. Since $x_+\notin\Sigma\x\{0\}$, we can take for large $n$, $\sigma_n\in\R$ the last $U$-entry time of $w_n$, i.e. $w_n(\sigma_n)\notin U$ and $w_n(s)\in U$ for $s>\sigma_n$. We note that $\sigma_n\to\infty$ as $n\to\infty$ and that the reparametized sequence $\sigma_n^*w_n$ is a gradient flow line of $\AA^\HH_{\sigma_n^*\varphi_{r_n} F}$ where $\sigma_n^*w_n(\cdot):=w_n(\cdot+\sigma_n)$ and $\sigma_n^*\varphi_{r_n}(\cdot):=\varphi_{r_n}(\cdot+\sigma_n)$. The new sequence $\sigma_n^*w_n$ also has a $C^\infty_{loc}$-convergent subsequence by Theorem \ref{thm:compactness of moduli} again and we denote by $z$ the limit gradient flow line. Since $r_n\to r$ and $\sigma_n\to\infty$, $\sigma_n^*\varphi_{r_n}$ $C^\infty_{loc}$-converges to the zero function, and thus $z$ is the gradient flow line of $\AA^\HH$.  Since $\sigma_n^*w_n\to z$ in $C^\infty_{loc}$-topology, we have
$$
E(z)=\int_{-\infty}^\infty||\p_s z||_m^2ds=\lim_{T\to\infty}\int_{-T}^T||\p_s z||_m^2ds\leq \lim_{T\to\infty}\limsup_{n\in\N}E(w_n)=\limsup_{n\in\N}E(w_n).
$$
We observe that $z(0)\notin U$ and the positive asymptotic end $z_+\in\Sigma\x\{0\}$ since $\Sigma\x\{0\}$ is a Morse-Bott component of $\Crit\AA^\HH$ (see \cite[Lemma 2.12]{AF1}) and hence $z$ is a non-constant gradient flow line of $\AA^\HH$. Thus the negative asymptotic end $z_-$ is a critical point of $\AA^\HH$; moreover it is not a trivial solution of \eqref{equation} since otherwise $z$ is a non-constant gradient flow line with zero energy $E(z)=0$. But this case is ruled out by the assumption that $||F||<\wp(\Sigma)$ as well. To be precise, with $z_-=(v,\eta)$, we can derive the following estimation which contradicts the definition of $\wp(\Sigma)$.
\beqn
0<|\Omega(v)|=|\AA^\HH_0(z_-)|=E(z)\leq\limsup_{n\in\N} E(w_n)\leq||F||<\wp(\Sigma).
\eeq

\noindent\underline{Case 2}. $r_n$ is unbounded.\\[-1.5ex]

Without loss of generality, we assume that $r_n\to\infty$ as $n\to\infty$. The limit of $\{w_n\}_{n\in\N}$ is a gradient flow line of $\AA^\HH_{F}$ since $\beta_\infty\equiv 1$. Then the asymptotic ends of the limit are critical points of $\AA^\HH_F$ which give rise to a leafwise intersection point of $\phi_F$. It contradicts our assumption and Case 2 is ruled out.\\[-2ex]

With $\sigma_n$ the first $U$-exit time of $w_n$, the case $x_-\notin\Sigma\x\{0\}$ is analogous. If $x_-=(q,0)\in\Sigma$ with $q\neq p$, as Case 1, there exists a gradient flow line of $\AA^\HH$ with asymptotic ends $(q,0)$ and $(p,0)$. But this cannot occur. Therefore we conclude that the moduli space $\M$ is compact. \\[-2ex]

Next, we regard the moduli space $\M$ as the zero set of a Fredholm section with index 1 of a Banach bundle over a Banach manifold as in \eqref{eq:Fredholm section}. Moreover, the Fredholm section is already transversal at the $(0,p,0)$ since $\Sigma$ is a Morse-Bott component by \cite[Lemma 2.12]{AF1}. Therefore we can perturb the Fredholm section away from $(0,p,0)$ (even if varying $J$, $(0,p,0)$ still solves the gradient flow equation) to obtain a transverse Fredholm section whose zero set is a compact one-dimensional smooth manifold with boundary $(0,p,0)$. But there is no one-dimensional manifold with a single boundary point. This finishes the proof of Claim 1. \hfill $\square$\\[-1ex]

\noindent\textbf{Step 2.} End of the proof of Theorem A.\\[1ex]
\noindent{Proof of Step 2}. In Step 2, our restricted contact coisotropic submanifold $\Sigma$ does not necessarily admit global coordinates. We consider a family of Hamiltonian tuples $\HH_\nu(t,x)=\chi(t)\GG_\nu(x)$, $\nu\in\N$ where $\HH_\nu=(H_{1,\nu},\dots,H_{k,\nu})$ and $\GG_\nu=(G_{1,\nu},\dots,G_{k,\nu})$ such that
\begin{enumerate}
\item $0<\epsilon_\nu<\min\{1/4k,\delta_0/2,\delta_1\}$ converges to zero as $\nu$ goes to infinity,
\item $G_{i,\nu}|_{U_{\delta_0}}=g_i(p_i)$ for some $g_i\in C^\infty(\R)$,
\item for $(x,\pp)\in\Sigma\x(-\delta_0,\delta_0)^k\cong U_{\delta_0}$,
\beq
G_{i,\nu}|_{U_{2\epsilon_\nu}-U_{\epsilon_\nu/2}}(x,\pp)= \left\{ \begin{array}{ll}
 p_i-\epsilon_\nu & \textrm{if}\,\, p_i>0\\[0.5ex]
 -p_i-\epsilon_\nu & \textrm{if}\,\, p_i<0, \end{array}\right.\;
\eeq
\item $G_{i,\nu}|_{M-U_{\delta_0}}=$ $constant$,
\item $\GG_\nu^{-1}(0)=\bigcup_{2^k}\Sigma\x(\pm\epsilon_\nu,\dots,\pm\epsilon_\nu)$.
\end{enumerate}
We note that
$$
X_{G_{i,\nu}}|_{\Sigma\x(\pm\epsilon_\nu,\dots,+\epsilon_\nu,\dots,\pm\epsilon_\nu)}=+X_{p_i},\quad
X_{G_{i,\nu}}|_{\Sigma\x(\pm\epsilon_\nu,\dots,-\epsilon_\nu,\dots,\pm\epsilon_\nu)}=-X_{p_i}.
$$
By construction, $\HH_\nu$ Poisson-commutes and Step 1 guarantees the existence of critical points $(v_\nu,\eta_\nu)$ lying on $\GG_\nu^{-1}(0)$ for sufficiently large $\nu$ because $||F||<\wp(\Sigma\x\{(\pm\epsilon_\nu,\dots,\pm\epsilon_\nu)\})$ for large $\nu\in\N$. For $(v_\nu,\eta_\nu)\in\Crit\AA^{\HH_\nu}_F$, $v_\nu$ lies on one of the components of $\GG_\nu^{-1}(0)$, say $v_\nu\subset\Sigma\x(\epsilon_\nu,\dots,\epsilon_\nu)$. According to Proposition \ref{prop:critical point answers question}, it holds that
$$
\phi_F^1\big(v_\nu(1/2)\big)=v_\nu(0)=\phi_{H_{1,\nu}}^{-\eta_{1,\nu}}\circ\cdots\circ\phi_{H_{k,\nu}}^{-\eta_{k,\nu}}\big(v_\nu(1/2)\big).
$$
Then the estimation \eqref{eq:step 1 in thm A} in Step 1 implies the following lemma.
\begin{Lemma}\label{lemma:uniform bound on time}
For $(v_\nu,\eta_\nu)\in\Crit\AA^{\HH_\nu}_F$, $\eta_{1,\nu},\dots,\eta_{k,\nu}$ are uniformly bounded in terms of $\lambda_1,\dots,\lambda_k$ and $F$.
\end{Lemma}
\begin{proof}
We estimate as in \eqref{eq:step 1 in thm A}: For all $i\in\{1,\dots,k\}$,
\bean
||F||&\geq\big|\AA^{\HH_\nu}_F(v_\nu,\eta_\nu)\big|\\
&=\Big|\int_0^1v^*\lambda_i+\int_0^1 \langle \eta,\HH_\nu\rangle(t,v_\nu(t))dt+\int_0^1F(t,v_\nu(t))dt\Big|\\
&=\Big|\int_0^1\lambda_i(v_\nu)\bigr(\sum_{j=1}^k\eta_{j,\nu} X_{H_{j,\nu}}(v_\nu)+X_F(t,v_\nu)\bigr)dt+\int_0^1F(t,v_\nu(t))dt\Big|\\
&=\frac{3}{4}|\eta_{i,\nu}|-\frac{1}{4(k-1)}\sum_{j\ne i}^k|\eta_{j,\nu}|-\Big|\int_0^1\lambda_i(v_\nu)\big(X_F(t,v_\nu)\big)+\int_0^1F(t,v_\nu(t))dt\Big|.\\
\eea
Therefore we conclude
\beqn
\frac{1}{2}\sum_{i=1}^k|\eta_{i,\nu}|\leq k\bigr(||F||+\max_{1\leq i\leq k}||\lambda_{i|{U_{\delta_0/2}}}||_{L^\infty}||X_F||_{L^\infty}+||F||_{L^\infty}\bigr).
\eeq
\end{proof}
The two sequences of points $\{v_\nu(0)\}_{\nu\in\N}$ and $\{v_\nu(1/2)\}_{\nu\in\N}$ converge up to taking a subsequence (still denoted by $v_\nu(0)$ and $v_\nu(1/2)$) and we denote by
\bean
x_0:=\lim_{\nu\to\infty}v_\nu(0),\quad x_{1/2}:=\lim_{\nu\to\infty}v_\nu\big(1/2\big).
\eea
Obviously $x_0$ and $x_{1/2}$ are points in $\Sigma$. Moreover we know that
\beq\label{eq:limit 1}
x_0=\lim_{\nu\to\infty}v_\nu(0)=\lim_{\nu\to\infty}\phi_F^1(v_\nu(1/2))=\phi_F^1(\lim_{\nu\to\infty}v_\nu(1/2))=\phi_F^1(x_{1/2}).
\eeq
Furthermore, due to Lemma \ref{lemma:uniform bound on time}, the limit  $\{\eta_{i,\nu}\}_{\nu\in\N}$ exists for all $i$, say
\beqn
\mathfrak{n}_i:=\lim_{\nu\to\infty}\eta_{i,\nu}.
\eeq
Thus we conclude that $x_0$ and $x_{1/2}$ lie on the same leaf:
\beq\label{eq:limit 2}
x_0=\lim_{\nu\to\infty}v_\nu(0)=\lim_{\nu\to\infty}\phi_{H_{1,\nu}}^{-\eta_{1,\nu}}\circ\cdots\circ\phi_{H_{k,\nu}}^{-\eta_{k,\nu}}(v_\nu(1/2))=
\phi_{H_{1}}^{-\mathfrak{n}_1}\circ\cdots\circ\phi_{H_{k}}^{-\mathfrak{n}_k}(x_{1/2}).
\eeq
It directly follows
\beqn
\phi_{H_{1}}^{-\mathfrak{n}_1}\circ\cdots\circ\phi_{H_{k}}^{-\mathfrak{n}_k}(x_{1/2})=\phi_F^1(x_{1/2})
\eeq
from \eqref{eq:limit 1} together with \eqref{eq:limit 2}. This completes the proof of Theorem A.
\hfill $\square$

\section{On contact coisotropic submanifolds}\label{unrestricted}
In this section, we explore the leafwise intersection problem on contact coisotropic submanifolds. As we mentioned in the introduction, a contact coisotropic submanifold notably differs from a restricted contact submanifold; for instance, an ambient symplectic manifold can be closed. For the leafwise intersection problem in the contact case, we shall again examine the perturbed Rabinowitz action functional $\AA^\HH_{F_s}$. However since the contact 1-forms are defined locally around $\Sigma$, we need a  constraint on the support of $X_F$. Once we achieve a uniform bound of Lagrange multipliers as before, then Theorem B follows from the proof of Theorem A. The main strategy is similar to \cite{Ka1}.\\[-1ex]

Let $(M,\om)$ be a symplectically aspherical symplectic manifold which is either closed or convex at infinity and $\Sigma$ be a closed contact coisotropic submanifold. We recall $\delta_0$ and $\FFF$:
\beqn
\delta_0:=\max\left\{
r\in\R\,\Bigg|\,\begin{aligned}&\textrm{ there exists a symplectic embedding }\\ \psi: U_r&=\{(q,\pp)\in\Sigma\x\R^k\,|\,|p_i|<r,\,\,  1\leq i\leq k\} {\hookrightarrow} M \end{aligned}
\;\;\right\}\eeq
Let $\psi_0:U_{\delta_0}\into M$ be a symplectic embedding and we identify $U_{\delta}$ with $\psi_0(U_{\delta})$ for all $0<\delta<\delta_0$.
\begin{Def}
A Hamiltonian function $F\in C^\infty(S^1\x M)$ is called {\em admissible} if $\Supp X_F\subsetneq U_{\delta_0}$. We abbreviate by $\FFF$ the set of all admissible Hamiltonian functions:
$$
\FFF:=\big\{F\in C_c^\infty(S^1\x M)\,|\,\Supp X_F\subsetneq U_{\delta_0}\big\}
$$
\end{Def}
We allow $s$-dependence on $F$ again as in subsection 4.2; in addition, $F_s\in\FFF$ for all $s\in\R$ and $\Supp X_{F_s}\subset\Supp X_{F_1}$. For $F_s\in\FFF$, we choose $\delta_2< \delta_0$ such that $\Supp X_{F_s}\subset U_{\delta_2}$. Then we modify the Hamiltonian functions $H_i(t,x)=\chi(t)p_i(x)$ to
\beqn
\widetilde H_i:M\rightarrow\mathbb{R}\,\,\ by \,\,\,
  \widetilde H_i = \left\{ \begin{array}{ll}
  H_i & \textrm {on  \,\,\,\,\,\,\,\,\,\,\,\,\,\,   $U_{{\delta_2}/2}$}\\
   constant & \textrm  {outside \,\,  $U_{\delta_2}.$ }\end{array}\right.\;
\eeq
By abuse of notation we write again $H_i$ for $\widetilde H_i$. In fact, the behavior of $H_i$ outside a small neighborhood of $\Sigma$ is not an important issue. Now we construct two cut-off functions to extend the contact 1-forms globally. We choose a function $\rho:\R\to\R$ satisfying
\begin{enumerate}\label{eq:cut-off}
\item $\rho(r)=r+\frac{3k+1}{4k},\quad r\in[-\delta_2,\delta_2],$
\item $\Supp\,\rho\subset(-\delta_0,\delta_0),$
\item $\rho'(r)\leq1+\varepsilon_0,\quad r\in\R$,\quad for some $\varepsilon_0>0$ satisfying
\beqn
\frac{(3k+1)/(4k)-\delta_2}{\delta_0-\delta_2}<1+\varepsilon_0;
\eeq
\end{enumerate}
and a function $\varrho:\R\to\R$ satisfying
\begin{enumerate}\label{eq:cut-off}
\item $\varrho(r)=r+\frac{1}{4k},\quad r\in[-\delta_2,\delta_2],$
\item $\Supp\,\varrho\subset(-\delta_0,\delta_0),$
\item $\varrho'(r)\leq1+\varepsilon_1,\quad r\in\R$,\quad for some $\varepsilon_1>0$ satisfying
\beqn
\frac{1/(4k)-\delta_2}{\delta_0-\delta_2}<1+\varepsilon_1.
\eeq
\end{enumerate}
The reason why we have such constraints on $\varepsilon_0$ and $\varepsilon_1$ is that our interval is not long enough. Then we extend the contact 1-forms $\alpha_i$'s to
\beqn
\beta_i(y):=\left\{\begin{array}{ll}
\rho(p_i)\alpha_i(y)+\sum_{j\neq i}^k\varrho(p_j)\alpha_j(y),\qquad & y=\psi(x,p_1,\dots,p_k)\in\Sigma\x(-\delta_0,\delta_0)^k  \\[2ex]
0,\qquad & y\in M\setminus(\Sigma\x(-\delta_0,\delta_0)^k)\end{array}\right.
\eeq
Recall that we have chosen a family of almost complex structures $\{J(s,t)\}_{(s,t)\in\R\x S^1}$ which splits on $U_{\delta_0}$ with respect to $TU_{\delta_0}=\xi\oplus\xi^\om$ as $J|_{\xi^\om}$ is an almost complex structure which interchanges Hamiltonian vector fields $X_{p_i}$ with $\frac{\partial}{\partial p_i}$ for $1\leq i\leq k$.

\begin{Prop}\label{prop:difference of bilinear forms}
For every $v\in TM$, the following inequality holds.
\beqn
d\beta_i(v,Jv)\leq(1+\varepsilon)\om(v,Jv)
\eeq
for $\varepsilon:=\max\{\varepsilon_0,\varepsilon_1\}$ and for all $i=1,\dots,k$.
\end{Prop}
\begin{proof}
Outside of $U_{\delta_0}$, the inequality is obvious since $d\beta$ vanishes but $\om(\cdot,J\cdot)=g(\cdot,\cdot)$ is positive definite. For $v\in TU_{\delta_0}$, we can write $v=v_1+v_2$ with respect to the decomposition $TU_{\delta_0}=\xi\oplus\xi^\om$. Since we have chosen $\rho$ and $\varrho$ so that $\rho'(r),\,\varrho'(r)\leq1+\varepsilon$
\bean
d\beta_i(v,Jv)&\leq\sum_{j\neq i}^k\varrho(p_j)\om_\Sigma(v_1,Jv_1)+\rho(p_i)\om_\Sigma(v_1,Jv_1)\\
&\quad+\sum_{j\neq i}^k\varrho'(p_j)dp_j\wedge\alpha_j(v_2,Jv_2)+\rho'(p_i)dp_i\wedge\alpha_i(v_2,Jv_2)\\
&\leq(1+\varepsilon)\Big(\bigr(\sum_{j=1}^k p_j+1\bigr)\om_\Sigma(v_1,Jv_1)+\sum_{j=1}^kdp_j\wedge\alpha_j(v_2,Jv_2)\Big)\\
&=(1+\varepsilon)\bigr(\om_\Sigma+\sum_{j=1}^kd(p_j\alpha_j)\bigr)(v,Jv)\\
&=(1+\varepsilon)\om(v,Jv).
\eea
\end{proof}
We define bilinear forms $\widehat m_i$'s on $T(\LLL\x\R^k)$ which are not necessarily positive definite.
\beqn
\widehat m_i\big((\hat v^1,\hat\eta^1),(\hat v^2,\hat\eta^2)\big):=\int_0^1d\beta_i(\hat v^1,J\hat v^2)dt+\hat\eta^1\hat\eta^2\;.
\eeq
We will use auxiliary action functionals
\beqn
\widehat\AA_{\beta_i,F_s}^{\HH}(v,\eta):=-\int_{D^2}\bar v^*d\beta_i-\int_0^1F_s(t,v(t))dt-\int_0^1\langle\eta,\HH\rangle(t,v(t))dt\;
\eeq
and
\beqn
\AA_i(v,\eta):=\widehat\AA_{\beta_i,F_s}^{\HH}(v,\eta)-\AA_{F_s}^{\HH}(v,\eta)=\int_{D^2}\bar v^*(\om-d\beta_i).\;
\eeq

\begin{Prop}\label{gradient w.r.t bilinear form}
For $(\hat{v},\hat{\eta}) \in T_{(v,\eta)}(\LLL \times \R^k)$, the following formula holds:
\beqn
 d\widehat\AA^{\HH}_{\beta_i,F_s}(v,\eta)[\hat v,\hat\eta] = \widehat m_i\Bigl(\nabla_m\AA^\HH_{F_s}(v,\eta),(\hat v,\hat\eta)\Bigr).\;
\eeq
\end{Prop}
\begin{proof}
On the region $\cl({U}_{\delta_2})$, we constructed $\rho$ and $\varrho$ so that
\beqn
d\beta_i=\Big(\sum_{j=1}^kp_j+1\Big)\om_\Sigma+\sum_{j=1}^k dp_j\wedge\alpha_j=\om_\Sigma+\sum_{j=1}^k d(p_j\alpha_j)=\om.
\eeq
On the other hand, we have chosen $\delta_2$ so that $X_{H_i}$'s and $X_{F_s}$ vanish outside of  $\cl({U}_{\delta_2})$. Therefore we have $i_{X_{H_i}}\om=i_{X_{H_i}}d\beta$ and $i_{X_{F_s}}\om=i_{X_{F_s}}d\beta$ and thus the following calculation proves the assertion.
\bea\label{eq:gradient wrt m and m hat}
d\widehat\AA^\HH_{\beta_i,F_s}(v,\eta)[\hat v,\hat\eta]
&=\int_0^1 d\beta_i(\p_t v,J\hat{v})-\om\bigr(\sum_{i=1}^k\eta_i{X_{H_i}(t,v)}+{X_{F_s}(t,v)},J\hat{v}\bigr) dt+\int_0^1\langle\hat{\eta},\HH(t,v)\rangle dt\\
&=\int_0^1 d\beta_i\bigr(\p_t v-\sum_{i=1}^k\eta_i{X_{H_i}(t,v)}-{X_{F_s}(t,v)},J\hat{v}\bigr) dt +\int_0^1\langle\hat{\eta},\HH(t,v)\rangle dt\\
&=\widehat m_i\Bigl(\nabla_m\AA^\HH_{F_s}(v,\eta),(\hat v,\hat\eta)\Bigr).
\eea
\end{proof}

\begin{Prop}\label{prop:bound on AA}
Along a gradient flow line $w=(u,\tau)$ of $\AA^\HH_{F_s}$,
\beqn
\big|\AA_i(w(s))\big| \leq \max\big\{\big|\AA_i(w_+)\big|,\big|\AA_i(w_-)\big|\big\}+\varepsilon E(w).\;
\eeq
where $\varepsilon=\max\{\varepsilon_0,\varepsilon_1\}$ and $w_\pm=\lim_{s\to\pm}w(s)$.
\end{Prop}
\begin{proof}
Using Proposition \ref{prop:difference of bilinear forms} and Proposition \ref{gradient w.r.t bilinear form},
\bea\label{eq:dA}
   \frac{d}{ds}\AA_i(w)
   =& d\widehat\AA^\HH_{\beta_i,F_s}(w)(\partial_s w) - d\AA^\HH_{F_s}(w)(\partial_s w)+\partial_s\widehat\AA_{F_s}^\HH(w)-\partial_s\AA_{F_s}^\HH(w)\\
   =& \int_0^1(\om-d\beta_i)(\partial_s u,J\partial_s u)dt+\int_0^1\partial_sF_sdt-\int_0^1\partial_sF_sdt\\
   \geq&-\int_0^1\varepsilon\om(\partial_s u,J\partial_s u)dt.
\eea
Integrating both sides of \eqref{eq:dA} with respect to $s$ from $-\infty$ to $s_0\in\R$, we obtain
\bean
\AA_i\big(w(s_0)\big)-\AA_i(w_-)
&= \int_{-\infty}^{s_0}\frac{d}{ds}\AA_i\big(w(s)\big)ds\\
&\geq-\varepsilon\int_{-\infty}^{s_0}\int_0^1\om\big(\partial_s u,J\partial_s u\big)dtds\\
&\geq-\varepsilon E(w).
\eea
In a similar way, we also have
\beqn
\AA_i\big(w(s_0)\big)-\AA_i(w_+)\leq\varepsilon E(w).
\eeq
The proposition follows from the above two inequalities.
\end{proof}

\begin{Prop}\label{prop:uniform bound of cutoff action functional}
$\widehat \AA_{\beta_i,F_s}^\HH$ is uniformly bounded along gradient flow lines of $\AA_{F_s}^\HH$.
\end{Prop}
\begin{proof}
By the definition of $\AA_i$, we know
\beqn
|\widehat\AA_{\beta_i,F_s}^\HH(w(s))|\leq|\AA_{F_s}^\HH(w(s))|+|\AA_i(w(s))|
\eeq
for a gradient flow line $w$ of $\AA^\HH_{F_s}$. But the righthand side of the above inequality is uniformly bounded by Proposition \ref{prop:uniform bound of action functional} and Proposition \ref{prop:bound on AA}.
\end{proof}
The next task is to find a uniform bound for the Lagrange multipliers. Using the fact that the action value of $\widehat\AA^\HH_{\beta_i,F_s}$ is uniformly bounded along gradient flow lines of $\AA^\HH_{F_s}$, we obtain a uniform bound of Lagrange multipliers as in the restricted contact case.

\begin{Lemma}\label{lemma:bound on eta when gradient is small again}
There exist $\epsilon>0$ and $C>0$ such that for $(v,\eta)\in\LLL\x\mathbb{R}^k$,
\beqn
||\nabla_m\AA^\HH_{F_s}(v,\eta)||_m<\epsilon \quad\textrm{  implies  }\quad |\eta_i|\leq C\Big(\sum_{j=1}^k|\widehat\AA^\HH_{\beta_j,F_s}(v,\eta)|+1\Big).
 \eeq
\end{Lemma}
\begin{proof}
We reformulate Step 1 in Lemma \ref{lemma:bound on eta when gradient is small} as below, then the same Steps 2 and 3 in Lemma \ref{lemma:bound on eta when gradient is small} finish the proof of the lemma.\\[0.5ex]
\textbf{Step 1:} Assume $v(t)\subset U_{\delta}$ for $t\in (1/2,1)$, where $\delta=\min\{1/4k,\delta_1,\delta_2/2\}$ (see \eqref{eq:near contact} for  $\delta_1$). Then there exists $C_0>0$ satisfying the following inequality.
\beqn
|\eta_i|\leq C_0\Big(\sum_{j=1}^k|\widehat\AA^\HH_{\beta_j,F_s}(v,\eta)|+||\nabla_m\AA^\HH_{F_s}(v,\eta)||_m+1\Big),\quad i=1,\dots,k.
\eeq
\noindent{Proof of Step 1.}
Since $H(t,x)=0$ for $t\in (0,1/2)$ and $v(t)\subset U_{\delta}$ for $t\in (1/2,1)$,
\bean
|\widehat\AA&^\HH_{\beta_i,F_s}(v,\eta)|
=\bigg|\int_0^1v^*\beta_i-\sum_{j=1}^k\eta_j\int_0^1H_j(t,v)dt-\int_0^1F_s(t,v)dt\bigg|\\
&\geq\bigg|\sum_{j=1}^k\eta_j\int_0^1\beta_i(v)\big(X_{H_j}(t,v)\big)dt\bigg|-\bigg|\int_0^1\beta_i(v)\big(X_{F_s}(t,v)\big)dt\bigg|-\bigg|\int_0^1F_s(t,v)dt\bigg|\\
&-\bigg|\int_0^1\beta_i(v)\big(\partial_t v-\sum_{j=1}^k\eta_jX_{H_j}(t,v)-X_{F_s}(t,v)\big)dt\bigg|-\bigg|\sum_{j=1}^k\eta_j\int_0^1H_j(t,v)dt\bigg|\\
&\geq\bigg|\eta_i\int_{1/2}^1\big(p_i+\frac{3k+1}{4k}\big)\alpha_i(v)\big(X_{H_i}(t,v)\big)\bigg|-\bigg|\sum_{j\ne i}^k\eta_j\int_{1/2}^1\big(p_j+\frac{1}{4k}\big)\alpha_j(v)\big(X_{H_j}(t,v)\big)\bigg|\\
&-\bigg|\int_0^1\beta_i(v)\big(X_{F_s}(t,v)\big)dt\bigg|-C_i||\nabla_m\AA^\HH_{F_s}(v,\eta)||_m-\bigg|\sum_{j=1}^k\eta_j\int_0^1H_j(t,v)dt\bigg|-\bigg|\int_0^1F_s(t,v)dt\bigg|\\
&\geq\frac{3}{4}|\eta_i|-\frac{2}{4k}\sum_{j\ne i}^k|\eta_j|-C_i||\nabla_m\AA^\HH_{F_s}(v,\eta)||_m-\delta\sum_{j=1}^k|\eta_j|-C_{i,F}
\eea
with $C_i:=||\beta_i|_{U_{\delta}}||_{L^\infty}<\infty$ and $C_{i,F}:=||F_s||_{L^\infty}+C_i||X_{F_s}||_{L^\infty}<\infty$. We have the above inequality for all $i=1,\dots,k$, thus we obtain
\beqn
\sum_{j=1}^k|\widehat\AA^\HH_{\beta_i,F_s}(v,\eta)|\geq\sum_{j=1}^k\Big(\frac{3}{4}-\frac{2(k-1)}{4k}-\frac{1}{4}\Big)|\eta_j|+\sum_{j=1}^k \bigr(C_j||\nabla_m\AA^\HH_{F_s}(v,\eta)||_m-C_{j,F}\bigr)
\eeq
and this proves Step 1. The same arguments as in Steps 2 and 3 in Lemma \ref{lemma:bound on eta when gradient is small} complete the proof of the lemma.
\end{proof}

\begin{Rmk}\label{Rmk:the reason we introduve the auxiliary Rabinowitz action functional}
The reason why we introduce the auxiliary action functionals is that we cannot achieve Step 1 in Lemma \ref{lemma:bound on eta when gradient is small} using only the perturbed Rabinowitz action functional. More precisely, since the one form $\alpha_i$ is not globally defined and $\bar v$ may go far away from $\Sigma$, we do not have the equality $\int_{D^2}\bar v^*\om=\int_0^1 v^*\alpha_i$.
\end{Rmk}

\begin{Lemma}\label{lemma:tau(sigma) again}
We have a bound on $o(\sigma,w,\epsilon)$ as follows:
\beqn
o(\sigma,w,\epsilon)\leq\frac{\AA^\HH_{F_s}(w_-)-\AA^\HH_{F_s}(w_+)+C_F}{\epsilon^2}.
\eeq
See \eqref{eq:def of tau and C_F} for the definitions of $o(\sigma,w,\epsilon)$ and $C_F$.
\end{Lemma}
\begin{proof}
It can be proved by exactly the same proof as Lemma \ref{lemma:tau(sigma)}.
\end{proof}

\begin{Thm}\label{thm:bound on eta again}
Assume that $w=(u,\tau)\in C^{\infty}(\mathbb{R},\LLL\x\mathbb{R}^k)$ is a gradient flow line of $\AA_{F_s}^\HH$ for which there exist $a \leq b$ such that
\beqn
\widehat\AA_{\beta_i,F_s}^\HH(w(s)),\,\,\AA_{F_s}^\HH(w(s))\in[a,b], \quad\textrm{for all }\, i=1,\dots,k,\,\, s\in \mathbb{R}.
\eeq
Then the $L^\infty$-norms of $\tau_i$'s are uniformly bounded.
\end{Thm}
\begin{proof}
Using Proposition \ref{prop:uniform bound of cutoff action functional}, Lemma \ref{lemma:bound on eta when gradient is small again} and Lemma \ref{lemma:tau(sigma) again}, we obtain
\bean
|\tau_i(\sigma)|&\leq|\tau_i(\sigma+o(\sigma,w,\epsilon))|+\int_\sigma^{\sigma+o(\sigma,w,\epsilon)}|\partial_s\tau_i(s)|ds\\
&\leq C\bigr(\sum_{j=1}^k|{\widehat\AA}^\HH_{\beta_j,F_s}(w)|+1\bigr)+o(\sigma,w,\epsilon)||H_i||_{L^\infty}\\
&\leq C(k\max\{|a|,|b|\}+1)+\Bigg(\frac{|b-a|+C_F}{\epsilon^2}\Bigg)||H_i||_{L^\infty}.
\eea
\end{proof}

\noindent\textbf{Proof of Theorem B.}
Since we got a uniform bound on Lagrange multipliers in Theorem \ref{thm:bound on eta again}, we again have Theorem \ref{thm:compactness of moduli} for $F_s\in\FFF$. Hence the proof of Theorem A proves Theorem B as well.
\hfill$\square$

\subsection{Proof of Corollary B}
Of course, $(S^1)^n$ in $((\C^*)^n,\om_\mathrm{std})$ is of restricted contact type; but we place this issue here since the technique used to prove Corollary B is similar to the technique used in the contact case.
To prove Corollary B, we need to check the compactness of gradient flow lines in the moduli space $\M$ in the proof of Theorem A. A uniform bound on the Lagrange multipliers can be established by the argument as in Section 4 since $(S^1)^n$ in $(\C^*)^n$ with the spherical symplectic form is of restricted contact type; here the spherical symplectic form is defined by $\om_\mathrm{std}=\sum_{i=1}^n dp_i\wedge d\theta_i$ where $p_i$ and $\theta_i$ are the coordinates on $(-1,\infty)$ and $S^1$ respectively. However, the problem is that $(C^*)^n$ is not convex at infinity. Thereby, we ought to show that gradient flow lines never escape to infinity.\\[-1.5ex]

We first fix a perturbation $F\in C^\infty_c(S^1\x (\C^*)^n)$ and note that  $\mathrm{cl}(\mathrm{Supp}X_F)$ is a compact subset of the $(\C^*)^n\cong (S^1)^n\x(-1,\infty)^n$. We recall the coordinate functions $p_i:\mathrm{cl}(\mathrm{Supp}X_F)\to(-1,\infty)$, and denote by
\beqn
\varrho^-:=\min_{i,x}p_i(x)\quad\& \quad\varrho^+:=\max_{i,x}p_i(x),\quad i\in\{1,\dots,k\},\;x\in(\C^*)^n.
\eeq
We choose a cut-off functions $\rho : \mathbb{R}\to\mathbb{R}$ which satisfies $\Supp\,\rho\subset[\varrho^--\epsilon,\varrho^++\epsilon]$ for any small $\epsilon>0$, $\rho(r)=r-|\varrho^-|$ on $[\varrho^-,\varrho^+]$, and $\rho'(r)\leq 1$ for all $r\in\mathbb{R}$. Then we have a global 1-form $\beta=\sum_{j=1}^n\rho(p_j)d\theta_j$. We modify the defining Hamiltonian functions $H_1,\dots,H_k$ to be constants outside of $[\varrho^-,\varrho^+]$. We again consider a one parameter family of perturbations $\{F_s\}_{s\in\R}$ such that $F_s$ varies only for $[-1,1]\subset\R$ and   $\Supp X_{F_s}\subset\Supp X_{F}$. We also define action functionals  $\widehat\AA_{F_s}^\HH$, $\AA$, and a bilinear form $\widehat m$ again:
\bean
\widehat\AA_{F_s}^{\HH}(v,\eta)&:=-\int_{D^2}\bar v^*d\beta-\int_0^1F_s(t,v(t))dt-\int_0^1\langle\eta,\HH\rangle(t,v(t))dt\;,\\
\AA(v,\eta)&:=\widehat\AA_{F_s}^{\HH}(v,\eta)-\AA_{F_s}^{\HH}(v,\eta)=\int_{D^2}\bar v^*(\om-d\beta),\\
\widehat m\big((\hat v^1,\hat\eta^1),(\hat v^2,\hat\eta^2)\big)&:=\int_0^1d\beta(\hat v^1,J\hat v^2)dt+\hat\eta^1\hat\eta^2\;,\quad (\hat v^1,\hat\eta^1),\,\,(\hat v^2,\hat\eta^2)\in T(\LLL\x\R^k).\\
\eea

\begin{Prop}\label{prop:two assertions}
For $(v,\eta) \in \LLL \times \R^n$ and  $(\hat{v},\hat{\eta}) \in
T_{(v,\eta)}(\LLL \times \R^n)$,
\begin{itemize}
 \item[(i)] $d\widehat\AA^\HH_{F_s}(v,\eta)(\hat v,\hat\eta) =
   \widehat m\Bigl(\nabla_m\AA^\HH_{F_s}(v,\eta),(\hat v,\hat\eta)\Bigr),$
\item[(ii)] $(m-\widehat m)\Bigl((\hat v,\hat\eta),(\hat v,\hat\eta)\Bigr)
   \geq 0.$
\end{itemize}
\end{Prop}

\begin{proof}
The proofs of the assertions (i) and (ii) are similar to the proofs of Proposition \ref{gradient w.r.t bilinear form} and Proposition \ref{prop:difference of bilinear forms}. For $v\in T((S^1)^n\x[\varrho^-,\varrho^+]^n)$, we write $v=v_1+v_2$ with respect to the decomposition $T((S^1)^n\x[\varrho^-,\varrho^+]^n)=\xi\oplus\xi^\om$. We estimate
\bean
d\beta(v,Jv)&=\sum_{i=1}^n\rho(p_i)d\alpha_i+\sum_{i=1}^n\rho'(p_i)dp_i\wedge\alpha_i\\
&\leq \max\Big\{\sum_{i=1}^n(p_i+\frac{1}{n})d\alpha_i(v_1,J_{\xi} v_1),0\Big\}+\om(v_2,J_{\xi^\om}v_2)\\
&\leq \max\Big\{\big(\om_\Sigma+\sum_{i=1}^nd(p_i\alpha_i)\big)(v_1,J_\xi v_1),0\Big\}+\om(v_2,J_{\xi^\om}v_2)\\
&\leq \om(v_1,J_\xi v_1)+\om(v_2,J_{\xi^\om}v_2)\\
&=\om(v,Jv).
\eea
Outside of the region $(S^1)^n\x[\varrho^-,\varrho^+]^n$, $d\beta$ vanishes but $\om$ is positive definite. This proves the assertion (ii).

In order to prove the assertion (i), it suffices to show that $i_{X_{F_s}}\om=i_{X_{F_s}}d\beta$ and $i_{X_{H_i}}\om=i_{X_{H_i}}d\beta$ for all $i=1,\dots,n$. We observe that $d\beta=\sum_{i=1}^n\rho'(p_i)dp_i\wedge d\theta_i$. Since $\rho'(p_i)=1$ on $[\varrho^-,\varrho^+]$, $\om=d\beta$. Outside of $[\varrho^-,\varrho^+]$, $X_{F_s}$ and $X_{H_i}$, $i=1,\dots,n$ vanish, thus $i_{X_{H_i}}\om=0=i_{X_{H_i}}d\beta$ and $i_{X_{F_s}}\om=0=i_{X_{F_s}}d\beta$ hold. The calculation \eqref{eq:gradient wrt m and m hat} proves the assertion (i).
\end{proof}

\begin{Cor}\label{cor:derivative of AA}
The functional $\AA$ is nondecreasing along gradient flow lines of $\AA_{F_s}^\HH$.
\end{Cor}
\begin{proof}
Using Proposition \ref{prop:two assertions}, we estimate for a gradient flow line $w$ of $\AA^\HH_{F_s}$,
\bean
\frac{d}{ds}\AA(w(s))
=& \frac{d}{ds}\bigg(\widehat\AA_{F_s}^\HH(w(s))\bigg)-\frac{d}{ds}\bigg(\AA_{F_s}^\HH(w(s))\bigg)\\
=& d\widehat\AA_{F_s}^\HH(w)(\partial_s w)+(\partial_s \widehat\AA_{F_s}^\HH)(w)-d\AA_{F_s}^\HH(w)(\partial_s w)-(\partial_s \AA_{F_s}^\HH)(w)\\
=& m\bigr(\nabla_m\AA_{F_s}^\HH(w),\nabla_m\AA_{F_s}^\HH(w)\bigr)-\widehat m\bigr(\nabla_m \AA_{F_s}^\HH(w),\nabla_m \AA_{F_s}^\HH(w)\bigr)\geq0.
\eea
\end{proof}

Now we consider the moduli space $\M$ and the family of perturbed  Rabinowitz action functionals $\AA^\HH_{\varphi_r(s)F}$ defined in the proof of Theorem A.

\begin{Cor}\label{cor:AA=0}
$\AA(w(s))\equiv0$ for all $(r,w)\in\M$.
\end{Cor}
\begin{proof}
We note that $\AA(w_+)=\AA(w_-)=0$ since the asymptotic ends $w_+$ and $w_-$ are constant solutions. Therefore the proof immediately follows from the previous corollary.
\end{proof}

\begin{Prop}
Assume $(r,w)=(r,u,\tau)\in\M$. Then $u\in C^\infty(S^1\x\R,(\C^*)^n)$ lies in $(S^1)^n \x[\varrho^-,\varrho^+]^n$.
\end{Prop}
\begin{proof}
Assume that $u(s,t)$ does not lie in $(S^1)^n\times[\varrho^--\epsilon,\varrho^++\epsilon]^n$ for $s_-<s<s_+$; this means that there exists a nonempty open subset $U\subset Z:=(s_-,s_+)\x S^1$ such that $u(s,t)\in (S^1)^n\x\big((-1,\varrho^--\epsilon)\cup(\varrho^++\epsilon,\infty)\big)^n$ for $(s,t)\in U$.\\
Using the previous corollary, we calculate
\bean
0=&\int_{s_-}^{s_+}\frac{d}{ds}\AA(w(s))\\
=&\int_{s_-}^{s_+}\int_0^1(\om-d\beta)(\partial_s u, J\partial_s u) dtds\\
=&\int_{Z-U}(\om-d\beta)(\partial_s u, J\partial_s u) dtds+\int_U\om(\partial_s u, J\partial_s u) dtds.
\eea
\\
But $(\om-d\beta)(\partial_s u, J\partial_s u)$ is bigger or equal to zero and $\int_U\om(\partial_s u, J\partial_s u) dtds>0$. Thus this case cannot occur and every gradient flow line of $\AA_{\varphi_r(s)F}^\HH$ in $\M$ lies in $(S^1)^n \times[\varrho^--\epsilon,\varrho^++\epsilon]^n$. Taking the limit $\epsilon\to 0$, this proves the proposition.
\end{proof}
\noindent\textbf{End of the proof of Corollary B.} Thanks to the previous proposition, we obtain a uniform $L^\infty$-bound on $u\in C^\infty(S^1\x\R,M)$ even though $(\C^*)^n$ is not convex at infinity. A uniform $L^\infty$-bound on the Lagrange multipliers follows from Theorem \ref{thm:bound on eta}; accordingly, a uniform $L^\infty$-bound on the derivatives of $u$ is also established, see the proof of Theorem \ref{thm:compactness of moduli}. Hence the proof of Theorem A goes through in this situation and we obtain a leafwise intersection point of $\phi_F\in\Ham_c((\C^*)^n,\om_{\mathrm{std}})$ for $||F||<\wp((S^1)^n)$. Moreover in this situation, a leafwise intersection point is nothing but a self intersection point and $\wp((S^1)^n)=\infty$ since every solution of \eqref{Bolle's eq} is non-contractible in $(\C^*)^n$. This finishes the proof of Corollary B.
\hfill $\square$\\

\section{On stable coisotropic submanifolds}
In this section we consider a stable coisotropic submanifold with global coordinates. At the end of the proof of Theorem D, the additional assumption on the existence of global coordinates will be removed. The existence of solution of \eqref{Bolle's eq} and the second inequality in  Theorem D were proved by Cieliebak-Frauenfelder-Paternain \cite{CFP} for separating stable hypersurfaces. Following their proof, we can extend their result to stable coisotropic submanifolds and simple observations slightly improve the theorem.\\[-1.5ex]

Let $\Sigma$ be a closed stable coisotropic submanifold with global coordinates in a symplectically aspherical symplectic manifold $(M,\om)$ which is either closed or convex at infinity. Suppose that $\Sigma$ is {\em displaced by} $F\in C^\infty_c(S^1\x M)$, i.e. $\phi_F(\Sigma)\cap\Sigma=\emptyset$. We consider again the smooth family of functions $\varphi_r\in C^\infty(\R,[0,1])$ defined in the proof of Theorem A. As before, we fix a point $p\in\Sigma$ and set
\beqn
\M:=\left\{(r,w)\in[0,\infty)\times C^{\infty}(\R,\LLL\times\R^k)\,\bigg|\, \begin{aligned}&\,\,\, w\text{ is a gradient flow line of $\AA^\HH_{\varphi_rF}$ with}\\&\!\lim_{s\to-\infty}w(s)=(p,0),\;\lim_{s\to\infty}w(s)\in\Sigma\x\{0\} \end{aligned}\right\}\;.
\eeq

\begin{Thm}\label{thm:compactness for stable case}
For $(r,w)\in\M$ where $w=(u,\tau)$, $\tau$ and $r$ are uniformly bounded.
\end{Thm}
 In the previous sections we showed how Rabinowitz Floer theory for hypersurfaces can be generalized to our set-up. Since the proof of Theorem \ref{thm:compactness for stable case} needs several technical lemmas and auxiliary action functionals as in the contact case, instead of giving the proof we refer to the reader \cite[Section 4.3]{CFP} or the earlier version of the present paper \cite[Section 6]{Ka4}.\\[-1.5ex]

\noindent\textbf{Proof of Theorem D.} As before, the proof proceeds in two steps;  we first prove the theorem when our stable coisotropic submanifold admits global coordinates and we remove this additional assumption in the next step. \\[1ex]
\noindent\textbf{Step 1.} We know that a sequence $\{(r_n,w_n)\}_{n\in\N}$ in $\M$ has a $C_{loc}^\infty$-convergent subsequence due to Theorem \ref{thm:compactness for stable case} together with the argument in the proof of Theorem \ref{thm:compactness of moduli}. We denote by $(r,w)$ the limit which is a gradient flow line of $\AA_{\varphi_rF}^\HH$. Again by compactness, $w$ asymptotically converges to $w_\pm=(v_\pm,\eta_\pm)\in\Crit\AA^\HH$ since $\varphi_r(\pm\infty)=0$. If $(r,w)\in\M$, the moduli space $\M$ is a one dimensional compact manifold with a single boundary point $\{(0,p,0)\}$ (after perturbing a Fredholm section as in the proof of Theorem A). However such a manifold does not exist and therefore one of the asymptotic ends $w_\pm$ of $w$ is a nontrivial solution of \eqref{Bolle's eq}. For simplicity, let us assume  $w_+\notin\Sigma\x\{0\}$. Following the notation from the proof of Theorem A, we consider $\sigma_n\in\R$ the last $U$-entry time. Then $\sigma_n^*w_n$ is a gradient flow line of $\AA^\HH_{\sigma_n^*\varphi_{r_n} F}$ and $C^\infty_{loc}$-converges to a non-constant gradient flow line $z$ of $\AA^\HH$ with $z(0)\notin U$ and $z_+\in\Sigma\x\{0\}$.\footnote{ Honestly speaking, we did not prove $C^\infty_{loc}$-convergence of $(r_n,\sigma_n^*w_n)$; but it follows from the proof of Theorem \ref{thm:compactness for stable case}.} By compactness and the energy estimate, $z_-=(v,\eta)\in\Crit\AA^\HH$ and $z_-$ is a nontrivial solution of \eqref{Bolle's eq}. Moreover, by \eqref{eq:uniform action bounds}, we have
\bean
-||F||\leq\AA^\HH_{\sigma_n^*\varphi_{r_n} F}(\sigma_n^*w_n(s))\leq||F||,\quad\forall s\in\R.
\eea
As $n$ goes to infinity, it holds that
\beq\label{eq:conclusion 1}
-||F||\leq\Omega(v)=\AA^\HH(z_-)\leq||F||
\eeq
for every Hamiltonian function $F\in C_c^\infty(S^1\x M)$ displacing $\Sigma$.
Since $\AA^\HH(z_+)=0$ and the action value of $\AA^\HH$ decreases along $z$,
\beq\label{eq:conclusion 2}
\big|\Omega(v)\big|=\big|\AA^\HH(z_-)\big|>0.
\eeq
\eqref{eq:conclusion 1} and \eqref{eq:conclusion 2} prove Theorem E provided that  $\Sigma$ admits global coordinates.\\[-1ex]

\noindent\textbf{Step 2.} Now we consider the situation that $\Sigma$ need not have global coordinates. We choose a family of Hamiltonian tuples $\HH_\nu(t,x)=\chi(t)\GG_{\nu}(x)$, $\nu\in\N$ where $\HH_\nu=(H_{1,\nu},\dots,H_{k,\nu})$ and $\GG_\nu=(G_{1,\nu},\dots,G_{k,\nu})$ such that
\begin{enumerate}
\item $0<\epsilon_\nu<\min\{1/4k,\delta_0/2,\delta_1\}$ converges to zero as $\nu$ goes to infinity,
\item $G_{i,\nu}|_{U_{\delta_0}}=g_i(p_i)$ for some $g_i\in C^\infty(\R)$,
\item for $(x,\pp)\in\Sigma\x(-\delta_0,\delta_0)^k\cong U_{\delta_0}$,
\beqn
G_{i,\nu}|_{U_{2\epsilon_\nu}-U_{\epsilon_\nu/2}}(x,\pp)= \left\{ \begin{array}{ll}
 p_i-\epsilon_\nu & \textrm{if}\,\, p_i>0\\[0.5ex]
 -p_i-\epsilon_\nu & \textrm{if}\,\, p_i<0, \end{array}\right.\;
\eeq
\item $G_{i,\nu}|_{M-U_{\delta_0}}=$ $constant$,
\item $\GG_\nu^{-1}(0)=\bigcup_{2^k}\Sigma\x (\pm\epsilon_\nu,\dots,\pm\epsilon_\nu)$.
\end{enumerate}
With this defining Hamiltonian tuple $\HH_\nu$, the argument in Step 1 still works and thus there exists $v_\epsilon\in\GG_\nu^{-1}(0)$ a solution of \eqref{Bolle's eq} satisfying $0<\Omega(v_\epsilon)\leq e(\GG_\nu^{-1}(0))$. Since $\GG_\nu^{-1}(0)$ is disconnected, $v_\epsilon$ lies in one of its connected components, say $v_\epsilon\subset\Sigma_\epsilon$. Since there is a diffeomorphism $\psi_\epsilon$ between $\Sigma_\epsilon$ and $\Sigma$,  $\psi_\epsilon(v_\epsilon)$ is a loop solving \eqref{Bolle's eq}, contractible in $M$ with $\Omega(\psi_\epsilon(v_\epsilon))=\Omega(v_\epsilon)>0$. Moreover if we have chosen sufficiently large $\nu$, $e(\Sigma)=e(\GG_\nu^{-1}(0))$. For simplicity, let us assume that $e(\Sigma)+\varepsilon<e(\GG_\nu^{-1}(0))$ for some small $\varepsilon>0$ and for all $\nu\in\N$; it means that there is $F\in C^\infty_c(S^1\x M)$ such that $||F||\in (e(\Sigma),e(\Sigma)+\varepsilon)$ such that $\phi_F(\Sigma)\cap\Sigma=\emptyset$; but if $\nu$ is big enough, $\phi_F$ also displaces $\GG_\nu^{-1}(0)$ and it contradicts $||F||< e(\GG_\nu^{-1}(0))$. Hence, we have proved that
$$
0<\Omega(\psi_\epsilon(v_\epsilon))=\Omega(v_\epsilon)\leq e(\GG_\nu^{-1}(0))=e(\Sigma).
$$
\hfill $\square$
\begin{Rmk}
If one succeeds in proving compactness of gradient flow lines of the perturbed Rabinowitz action functional in the stable case, Theorem D is an immediate consequence of the invariance property of Rabinowitz Floer homology.
\end{Rmk}

\section{Rabinowitz Floer Homology}\label{Rabinowitz Floer Homology}

For hypersurfaces, \cite{CFP,AF1} proved that the (perturbed) Rabinowitz action functional is generically Morse-Bott (Morse). Their argument undeniably continuous to hold in our set-up. That is, $\AA_F^\HH$ is Morse for a generic perturbation $F\in C_c^\infty(S^1\x M)$. In the restricted contact case, furthermore, we know that gradient flow lines of the Rabinowitz action functional are compact modulo breaking (see (F1) and (F2) below) due to Theorem \ref{thm:compactness of moduli}. Thus we can define Floer homology of $\AA_F^\HH$ as usual. We denote this homology by  $\HF(\AA^\HH_F)$. If there is no perturbation, i.e. $F\equiv0$, $\AA^\HH$ is never Morse since there is a $S^1$-symmetry coming from time-shift on the critical points set. However $\AA^\HH$ is Morse-Bott for a generic coisotropic submanifold \footnote{ Since Rabinowitz Floer homology is invariant under homotopies there is no loss of generality in assuming $\AA^\HH$ is Morse-Bott, see \cite{CFP}.}, thus we can define Morse-Bott homology of $\AA^\HH$ by counting gradient flow lines with cascades, see \cite{F}. We define Rabinowitz Floer homology by $\RFH(\Sigma,M)=\HF(\AA^\HH)$. As one expects, these two Floer homologies, $\HF(\AA^\HH_F)$ and $\RFH(\Sigma,M)$, are isomorphic by the standard continuation argument in Floer theory. We only treat the closed restricted contact coisotropic submanifold $\Sigma$ in this section and refer to Remark \ref{rmk:RFH for contact and stable} for the other cases. As before, $(M,\om)$ is an exact symplectic manifold being convex at infinity with a family of $\om$-compatible almost complex structures $J=J(s,t)$.

\subsection{Boundary Operator}
We can assign some index to critical points of $\AA^\HH_F$, namely the transverse Conley-Zehnder index.\footnote{ We still can define Floer homology of $\AA_F^\HH$ without this index.} But we omit the definition, referring the reader to \cite{BO2,CF,MP}. We denote the index by
$$
\mu:\Crit\AA^\HH_F\pf \Z.
$$
Here we assumed that the first Chern class $c_1$ vanishes over $\pi_2(M)$ for simplicity; otherwise the index $\mu$ is well-defined modulo $2N$ where $N$ is the minimal Chern number of $(M,\om)$.

Let $\M_J(w_-,w_+)$ be the moduli space of gradient flow lines of $\AA_{F}^\HH$ with asymptotic ends $w_\pm\in\Crit\AA^\HH_{F}$.
$$
\M_J(w_-,w_+):=\left\{(u,\tau)\in C^\infty(\R\x S^1,M)\x C^\infty(\R,\R^k)\,\Bigg|\,\begin{aligned} (u,\tau) \textrm{ solves } \eqref {eq:Floer's interpretation},\,\\ \lim_{s\to\pm\infty}(u,\tau)=w_\pm\;\end{aligned}\right\}.
$$
In order to show that $\M_J(w_-,w_+)$ is a finite dimensional smooth manifold, we interpret it as the zero set of a Fredholm section of a Banach bundle over a Banach space. Let $\PP(w_-,w_+)$ be the Banach manifold given by
$$
\PP(w_-,w_+):=\big\{(u,\tau)\in W^{1,2}(\R\x S^1, M)\x W^{1,2}(\R,\R^k)\,\big|\, \lim_{s\to\pm\infty}(u,\tau)=w_\pm\big\}
$$
and $\EE$ be the Banach bundle over $\PP(w_-,w_+)$ whose fibre at $(u,\tau)\in\PP(w_-,w_+)$ is
$$
\EE_{(u,\tau)}:=L^2(\R\x S^1, u^*TM\x \tau^*T\R^k).
$$
Then the moduli space $\M(w_-,w_+)$ is the zero set of the section
\beq\label{eq:Fredholm section}
s_J:\PP(w_-,w_+)\pf\EE,\quad s_J(u,\tau)=\big(\bar\p_{\HH,F,J}(u),\bar\p_1(\tau_1),\cdots,\bar\p_k(\tau_k)\big)
\eeq
defined by
\bean\left.\begin{aligned}
\bar\p_{\HH,F,J}(u)&=\big(\p_s u+J(s,t,u)\p_tu-\sum_{i=1}^k\eta_iX_{H_i}(t,u)-X_{F_s}(t,u)\big)\\
\bar\p_i(\tau_i)&=\p_s\tau_i-\int_0^1 H_i(t,u)dt,\qquad 1\leq i\leq k
\end{aligned}\;\;\right\}\eea
where $\tau=(\tau_1,\dots,\tau_k)$. It turns out that this section is Fredholm. Then we regard the moduli space as the zero set of this section, $\M_J(w_-,w_+)=s^{-1}_J(0)$. Let
$$
Ds_J(u,\tau):T_{(u,\tau)}\PP(w_-,w_+)\pf \EE_{(u,\tau)}
$$
be the vertical differential of $s_J$ at $(u,\tau)$. It is known that $Ds_J(u,\tau)$ is surjective for generic $\om$-compatible almost complex structures $J$ and for any $(u,\tau)\in s_J^{-1}(0)$, see \cite[Section 5]{FHS} and \cite{BO1}. This transversality issues (surjectivity of $Ds_J(u,\tau)$) can now also be settled using the frame work of polyfolds developed by Hofer-Wysocki-Zehnder \cite{HWZ1,HWZ2,HWZ3}. Thus we perturb the section $s_J$ (varying $J$ slightly) so that $Ds_J(u,\tau)$ is surjective and the implicit function theorem yields that $s^{-1}_J(0)=\M_J(w_-,w_+)$ is a smooth finite dimensional manifold. Moreover the dimension of the moduli space $\M_J(w_-,w_+)$ coincides with the dimension of the kernel of $Ds_J$ which in turn is same as the Fredholm index of $s_J$ since it is surjective; besides, the Fredholm index of $s_J$ can be computed in terms of the indices of $\mu(w_-)$ and $\mu(w_+)$ using the spectral flow \cite{RS,BO2,CF}. In conclusion, we have the identity
$$
\dim\M_J(w_-,w_+)=\mu(w_-)-\mu(w_+),\quad w_\pm\in\Crit\AA^\HH_F.
$$
We suppress the subindex $J$ in $\M_J(w_-,w_+)$ for notational convenience. We divide out the $\R$-action on $\M(w_-,w_+)$ defined by shifting the gradient flow lines in the $s$-variable. Then we obtain the moduli space of unparametrized gradient flow lines which we denote by
\beqn
\widehat\M(w_-,w_+):=\M(w_-,w_+)/\R.
\eeq

For the compactification of the moduli space $\M(w_-,w_+)$, we recall the {\em Floer-Gromov convergence}: A sequence $\{(u^\nu,\tau^\nu)\}_{\nu\in\mathbb{N}}$ in $\M(w_-,w_+)$ is said to Floer-Gromov converge to a broken flow line $\{(u_j,\tau_j)\}_{j=1}^m$ where
$z_0,\dots,z_m\in \Crit\AA^\HH_{F_s}$ with $z_0=w_-$ and $z_m=w_+$ and
$$
(u_j,\tau_j)\in\M(z_{j-1},z_j),\quad j\in\{1,\dots, m\}.
$$
if there exist $\sigma_j^\nu\in\R$ such that reparametrized sequences $(u^\nu,\tau^\nu)(\sigma^\nu_j+\cdot)$ converge to $(u_j,\tau_j)$ for all $j\in\{1,\dots, m\}$ in the $C^\infty_{loc}$-topology. The following statements are the key ingredients for boundary operators of various Floer homologies, including Rabinowitz Floer homology.

\begin{itemize}
\item[(F1)] The moduli space $\M(w_-,w_+)$ is a one dimensional compact smooth manifold with respect to the topology of Floer-Gromov convergence when $\mu(w_-)-\mu(w_+)=1$.\footnote{ Without help of the index, we can rephrase that the one dimensional component of $\M(w_-,w_+)$ is a compact smooth manifold.} Accordingly, $\widehat\M(w_-,w_+)$ is a finite set.
\item[(F2)] Let ${\widehat\M_c(w_-,w_+)}$ be the compactification of $\widehat\M(w_-,w_+)$ with respect to the topology of Floer-Gromov convergence. If $\mu(w_-)-\mu(w_+)=2$, ${\widehat\M_c(w_-,w_+)}$ is a compact one-dimensional manifold whose boundary is
\beq\label{eq:broken flow lines}
\p{\widehat\M_c(w_-,w_+)}=\bigcup_z\widehat\M(w_-,z)\x\widehat\M(z,w_+)
\eeq
where the union runs over $z\in\Crit\AA^\HH_F$ with $\mu(w_-)-\mu(z)=1$.
\end{itemize}
(i) follows from the elliptic bootstrapping argument as discussed in Theorem \ref{thm:compactness of moduli}, see also Floer's beautiful paper \cite{Fl2}. (ii) is proved by Floer's gluing theorem \cite{Fl1}.

We denote by $\Crit_q(\AA^\HH_F)$ the set of critical point of $\AA^\HH_F$ of index $q\in\Z$, i.e. $\mu((v,\eta))=q$ for $(v,\eta)\in\Crit_q(\AA^\HH_F)$. We define a $\Z/2$-vector space
\beqn
\CF_q(\AA^\HH_F):=\Big\{\xi=\!\!\!\!\!\!\!\sum_{(v,\eta)\in\Crit_q\AA^\HH_F}\!\!\!\xi_{(v,\eta)} (v,\eta)\,\Big|\,\xi_{(v,\eta)}\in\Z/2\Big\}
\eeq
where $\xi_{(v,\eta)}$ satisfies the finiteness condition:
\beqn
\#\bigr\{{(v,\eta)}\in\Crit_q\AA^\HH_F\,\bigr|\,\xi_{(v,\eta)}\ne0,\,\,\AA^\HH_F(v,\eta)\geq\kappa\bigr\}<\infty,\quad\forall\kappa\in\R.
\eeq
We denote by $n(w_-,w_+)$ be the parity of elements of the finite set  $\widehat\M(w_-,w_+)$ when $\mu(w_-)-\mu(w_+)=1$, see (F1) above. Then the boundary operators $\{\p_q\}_{\{q\in\Z\}}$ are defined by
\bean
\p_q:\CF_q(\AA_F^\HH)&\pf\CF_{q-1}(\AA_F^\HH)\\
w_-\in\Crit_q\AA^\HH_F&\longmapsto\!\!\!\!\!\!\!\sum_{w_+\in\Crit_{q-1}\AA^\HH_F}\!\!\!\!\! n(w_-,w_+)\cdot w_+.
\eea
Due to (F2), we know $\p_{q-1}\circ\p_q=0$ (in $\Z/2$) so that $(\CF_*(\AA_F^\HH),\p_*)$ is a chain complex indeed. We define Rabinowitz Floer homology by
\beqn
\HF_q(\AA^\HH_F):=\H_q(\CF_*(\AA_F^\HH),\p_*),\quad\RFH_q(\Sigma,M):=\HF_q(\AA^\HH).
\eeq
To be exact, since $\AA^\HH$ is Morse-Bott, $\HF(\AA^\HH)$ is defined by Frauenfelder's Morse-Bott homology \cite[Appendix A]{F}. We note that $\Crit\AA^\HH$ consists of $\Sigma$ and circles. We pick a Morse function $f$ on $\Crit\AA^\HH$ and then the boundary operator for $\HF(\AA^\HH)$ is defined by counting gradient flow lines of $\AA^\HH$ (cascades) with gradient flow lines of $f$. One thing we use here is that if there is no nontrivial solution of \eqref{Bolle's eq}, $\Crit\AA^\HH=\Sigma$ and thus there are no cascades since the energy of each cascade would be positive. Thus if this is the case, $\HF(\AA^\HH)\cong \H(\Sigma)$.

\subsection{Continuation Homomorphism}
Given any two time-dependent Hamiltonian functions $F$ and $K$, we consider the homotopies $D_s^\pm\in C^\infty(S^1\x M)$, $s\in\R$,
\beqn
D_s^+(t,x):=K(t,x)+\varphi_+(s)\big(F(t,x)-K(t,x)\big)
\eeq
and
\beqn
D_s^-(t,x):=K(t,x)+\varphi_-(s)\big(F(t,x)-K(t,x)\big)
\eeq
where $\varphi_\pm\in C^\infty(\R,[0,1])$ are cut-off functions defined by
\beqn
\varphi_+(s)=\left\{\begin{aligned}0\quad & s\leq -1\\
1\quad & s\geq 1\end{aligned}\right.\qquad
\varphi_-(s)=\left\{\begin{aligned}1\quad & s\leq -1\\
0\quad & s\geq 1\end{aligned}\right.
\eeq
Now we consider the time-dependent version of the gradient flow equation:
\beq\label{eqn:gradient flow equation for G_s}\left.
\begin{aligned}
&\partial_su+J_s(t,u)\bigr(\partial_tu-\sum_{i=1}^k \tau_i X_{H_i}(t,u)-X_{D_s^+}(t,u)\bigr)=0\\
&\partial_s\tau_i-\int_0^1H_i(t,u)dt=0, \qquad 1\leq i\leq k.
\end{aligned}
\;\;\right\}
\eeq
The solutions of \eqref{eqn:gradient flow equation for G_s} with the asymptotic conditions form the following moduli space:
\beqn
\M(w_K,w_F)\!:=\!\left\{w\in C^\infty(\R\x S^1,M)\!\x\! C^\infty(\R,\R^k)\,\Bigg|\begin{aligned}\,w=(u,\tau) \textrm{ solves } \eqref{eqn:gradient flow equation for G_s}\textrm{ with }\quad  \\ \lim_{s\to\pm\infty}\!\!w(s)=w_{F/K}\in\Crit\AA^\HH_{F/K}\, \end{aligned}\right\}.
\eeq
As we discussed in the previous subsection, it is also a well-known fact in Floer theory that the moduli space $\M(w_K,w_F)$ is a smooth manifold of dimension $\mu(w_K)-\mu(w_F)$ for a generic homotopy. In particular, it is known that $\M(w_K,w_F)$ is a finite set when $w_K$ and $w_F$ have the same index and thus we denote the parity of $\M(w_K,w_F)$ by $n(w_K,w_F)$ if this is the case. Then we define the continuation homomorphism as follows.
\bean
\Phi_K^F:\CF_q(\AA_K^\HH)&\pf\CF_q(\AA_F^\HH)\\
w_K\in\Crit_q\AA^\HH_K&\longmapsto\!\!\!\!\!\!\!\sum_{w_F\in\Crit_q\AA^\HH_F}\!\!\!\!\! n(w_K,w_F)\cdot w_F.
\eea
In the same way, we also define $\Phi_F^K:\CF(\AA_F^\HH)\to\CF(\AA_K^\HH)$ using the other homotopy $D_s^-$. Then we obtain the invariance property of Rabinowitz Floer homology via the continuation homomorphisms using a  homotopy of homotopies $D^r_s(t,x):=K(t,s)+\varphi_r(s)(F(t,x)-K(t,x))$ where $\varphi_r:\R\to[0,1]$, $r\in\R$ and $\varphi_r=\varphi_\pm$ if $\pm r\geq 1$, see \cite[Section 3.4]{S} \footnote{ Here we again make use of Floer-Gromov compactness and Floer's gluing theorem.}:
\begin{Thm}\label{Thm:invariance of RFH}
Rabinowitz Floer homology is independent of the choice of perturbations up to canonical isomorphism. In particular, it holds that
\beqn
\RFH(\Sigma,M)\cong\HF(\AA^\HH_F),\quad F\in C_c^\infty(S^1\x M).
\eeq
\end{Thm}

\noindent\textbf{Proof of Theorem E.}
Suppose that there are no leafwise intersection points for some $\phi_F\in\Ham_c(M,\om)$. Then the set $\Crit\AA_F^\HH$ is empty since otherwise a critical point of $\AA_F^\HH$ gives rise to a leafwise intersection point. Thus $\HF(\AA^\HH_F)=0$ and Theorem  \ref{Thm:invariance of RFH} proves (i).

If there are only trivial solutions of \eqref{Bolle's eq}, no cascades appear in the boundary operator of Morse-Bott homology. Thus  Rabinowitz Floer homology is isomorphic to Morse homology of $\Sigma$ and hence to singular homology of $\Sigma$. This proves (iii).

Suppose there are only trivial solutions of \eqref{Bolle's eq}. Due to (iii), we know that Rabinowitz Floer homology is isomorphic to singular homology of $\Sigma$. Singular homology of $\Sigma$ never vanishes, but Rabinowitz Floer homology vanishes by (i) since $\Sigma$ is displaceable. This contradiction proves (ii).
\hfill $\square$\\[-0.5ex]

For the later purpose, we compare the action values of $\AA^\HH_{K}$ and $\AA^\HH_{F}$:
\begin{Prop}\label{Prop:compare action values}
If the moduli space $\M(w_K,w_F)$ is not empty,
\beqn
\AA^\HH_F(w_F)\leq\AA^\HH_K(w_K)+||F-K||_-.
\eeq
\end{Prop}
\begin{proof}
We pick $w\in\M(w_K,w_F)$ and estimate its energy:
\bean
0&\leq E(w)\\
&=-\int_{-\infty}^{\infty}d\AA^\HH_{D^+_s}(w(s))[\p_s w]ds\\
&=-\int_{-\infty}^{\infty}\frac{d}{ds}\bigr(\AA_{D^+_s}^\HH(w(s))\bigr)ds-
\int_{-\infty}^{\infty}\int_0^1{\varphi_+}'(s)\big(F(t,w(s))-K(t,w(s))\big)dtds.\\
&\leq\AA^\HH_{D^+_{-\infty}}(w_K)-\AA^\HH_{D^+_\infty}(w_F)
-\int_{-\infty}^{\infty}{\varphi_+}'(s)\int_0^1\big(F(t,w(s))-K(t,w(s))\big)dtds.\\
&\leq\AA^\HH_K(w_K)-\AA^\HH_F(w_F)+||F-K||_-.
\eea
\end{proof}

\subsection{Filtered Rabinowitz Floer Homology}
For $a<b\in\R$ which are not critical values of $\AA^\HH_F$, we define the $\Z/2$-vector space
\beqn
\CF_q^{(a,b)}(\AA^\HH_F):=\Crit_q^{(a,b)}(\AA^\HH_F)\otimes\Z/2
\eeq
where
\beqn
\Crit_q^{(a,b)}(\AA^\HH_F):=\big\{(v,\eta)\in\Crit_q\AA^\HH_F\,\big|\,\AA^\HH_F(v,\eta)\in(a,b)\big\}.
\eeq
Then $\big(\CF_*^{(-\infty,b)}(\AA^\HH_F),\p_*^b\big)$ is a sub-complex of $\big(\CF_*(\AA^\HH_F),\p_*\big)$ since gradient flow lines of $\AA^\HH_F$ flow downhill. Here $\p_*^b:=\p_*|_{\CF_*^{(-\infty,b)}}$. There are canonical homomorphisms
\beqn
i_a^{b,c}:\CF_q^{(a,b)}(\AA^\HH_F)\pf\CF_q^{(a,c)}(\AA^\HH_F),\qquad a\leq b\leq c
\eeq
and
\beqn
\pi^c_{a,b}:\CF_q^{(a,c)}(\AA^\HH_F)\pf\CF_q^{(b,c)}(\AA^\HH_F),\qquad a\leq b\leq c.
\eeq
$i_a^{b,c}$ is an inclusion and $\pi^c_{a,b}$ is a projection along $\CF_q^{(a,b)}(\AA^\HH_F)$. We note that
\beqn
\CF_q^{(a,c)}(\AA^\HH_F)=\CF_q^{(a,b)}(\AA^\HH_F)\oplus\CF_q^{(b,c)}(\AA^\HH_F),
\eeq
We suppress the indices $a$, $b$, and $c$ if there is no confusion. The short exact sequence
\beqn
0\pf\CF_q^{(-\infty,a)}(\AA^\HH_F)\stackrel{i}{\pf}\CF_q^{(-\infty,b)}(\AA^\HH_F)\stackrel{\pi}{\pf}\CF_q^{(a,b)}(\AA^\HH_F)\pf0,
\eeq
gives rise to a boundary operator $\p_{a*}^b$ on $\CF_*^{(a,b)}(\AA^\HH_F)$ and this induces a homology group, namely
{\em filtered Rabinowitz Floer homology}:
\beqn
\HF_q^{(a,b)}(\AA^\HH_F)=\H_q(\CF_*^{(a,b)}(\AA^\HH_F),\p_{a*}^b).
\eeq
More generally for $a\leq b\leq c$, we have
\beqn
0\pf\CF_q^{(a,b)}(\AA^\HH_F)\stackrel{i}{\pf}\CF_q^{(a,c)}(\AA^\HH_F)\stackrel{\pi}{\pf}\CF_q^{(b,c)}(\AA^\HH_F)\pf0.
\eeq
The canonical homomorphisms $i$, $\pi$, and the boundary map $\p$ are compatible with each other so that they induce canonical homomorphisms on homology level. Thus we have
\beqn
\cdots\stackrel{\delta}{\pf}\HF^{(a,b)}_q(\AA^\HH_F)\stackrel{i_*}{\pf}\HF^{(a,c)}_q(\AA^\HH_F)\stackrel{\pi_*}{\pf}
\HF^{(b,c)}_q(\AA^\HH_F)\stackrel{\delta}{\pf}\HF^{(a,b)}_{q-1}(\AA^\HH_F)\stackrel{i_*}{\pf}\cdots.
\eeq
where $\delta$ is the connecting homomorphism.

\begin{Cor}\label{Cor:canonical homomorphism for the filtered case}
The canonical homomorphism for the filtered case is defined by
\beqn
(\Phi_K^F)_*:\HF_q^{(a,b)}(\AA^\HH_K)\pf\HF_q^{(a-||F-K||_-,b+||F-K||_-)}(\AA^\HH_F)
\eeq
\end{Cor}
\begin{proof}
This is a well-known fact in Floer theory; it follows from the comparison of the action values of $\AA^\HH_K$ and $\AA^\HH_F$, see Proposition \ref{Prop:compare action values}.
\end{proof}

\subsection{Local Rabinowitz Floer Homology}
All of the lemmas and the propositions in this subsection were established for hypersurfaces in \cite{AF1}. Without doubt, their arguments continue to hold in our situation, but for the sake of completeness we outline the arguments.

For $||F||<\wp(\Sigma)$, we define
\beqn
\Crit_\mathrm{loc}(\AA^\HH_F):=\Big\{(v,\eta)\in\Crit\AA^\HH_F\,\Big|\,-||F||_+\leq\AA^\HH_F(v,\eta)\leq||F||_-\Big\}\;.
\eeq
We note that the set $\Crit_\mathrm{loc}(\AA^\HH_F)$ is finite. This follows from the Arzela-Ascoli theorem since the Lagrange multipliers $\eta_i$'s are uniformly bounded according to Theorem \ref{thm:bound on eta}.  We define the finite dimensional $\Z/2$ vector space
\beqn
\CF_\mathrm{loc}(\AA^\HH_F):=\Crit_\mathrm{loc}(\AA^\HH_F)\otimes\Z/2\;.
\eeq
$(\CF_\mathrm{loc}(\AA^\HH_F),\partial_\mathrm{loc})$ is a chain complex  and the {\em local Rabinowitz Floer homology} is defined by
\beqn
\HF_\mathrm{loc}(\AA^\HH_F):=\H(\CF_\mathrm{loc}(\AA^\HH_F),\partial_\mathrm{loc}).
\eeq

\begin{Prop}\label{Prop:betti bumber inequalities}
If $||F||<\wp(\Sigma)$, the following inequalities hold.
\beqn
\nu_{\mathrm{leaf}}(\phi_H)\geq\dim\CF_\mathrm{loc}(\AA_F^\HH)\geq\dim\HF_\mathrm{loc}(\AA^\HH_F)\;.
\eeq
Here $\nu_{\mathrm{leaf}}(\phi_F)$ is the number of leafwise intersection points of $\phi_F\in\Ham_c(M,\om)$.
\end{Prop}
\begin{proof}
We briefly sketch the proof and refer to \cite[Lemma 2.19]{AF1} for details. The last inequality is obvious. For the first inequality, it suffices to show that different critical points of $\AA^\HH_F$ give rise to different leafwise intersection points. If two critical points $(v,\eta)\ne(v',\eta')\in\Crit_\mathrm{loc}\AA^\HH_F$ give rise to the same leafwise intersection point, then $\gamma:=\underline v'|_{[1/2,1]}\# v|_{[1/2,1]}$, where $\underline v(t)=v(1-t)$ and $\#$ is the path catenation operator, is a closed orbit solving \eqref{Bolle's eq}. Moreover a close look at $\gamma$ reveals that $\Omega(\gamma)\leq||F||<\wp(\Sigma)$ which contradicts the definition of $\wp(\Sigma)$.
\end{proof}

\begin{Prop}\label{Prop:RFHloc=H}
Local Rabinowitz Floer homology of $\AA^\HH$ is isomorphic to singular homology of $\Sigma$, i.e.
\beqn
\H(\Sigma;\Z/2)\stackrel{\Theta}{\cong}\HF_\mathrm{loc}(\AA^\HH)\;.
\eeq
\end{Prop}
\begin{proof}
The set $\Crit_\mathrm{loc}\AA^\HH$ consists of critical points of $\AA^\HH$ whose action values are zero. According to Proposition \ref{Lem:eta's are same} and \eqref{eq:eta is period}, the action value of a critical point is equal to the Lagrange multiplier and thus $\Crit_\mathrm{loc}\AA^\HH=\Sigma$. Therefore no cascades appear in the boundary operator and $\HF_\mathrm{loc}(\AA^\HH)$ is isomorphic to Morse homology of $\Sigma$.
\end{proof}

The lemma below directly follows from the definition of $\wp(\Sigma)$.
\begin{Lemma}\label{Lem:filtered RHM rel to local RFH}
For any $(a,b)\subset(-\wp(\Sigma),\wp(\Sigma))$, we have an isomorphism
\beqn
\HF^{(a,b)}(\AA^\HH)\cong\HF_\mathrm{loc}(\AA^\HH).
\eeq
\end{Lemma}

\begin{Prop}\label{Prop:injective homomorphism}
If $||F||<\wp(\Sigma)$, there exists an injective homomorphism
\beqn
\iota:\H(\Sigma;\Z/2)\pf\HF_\mathrm{loc}(\AA^\HH_F)\;.
\eeq
In particular, $\dim\HF_\mathrm{loc}(\AA^\HH_F)\geq\dim\H(\Sigma;\Z/2)$.
\end{Prop}
\begin{proof}
We pick $a\in\R$ with $0<a<||F||<\wp(\Sigma)$ then using the continuation homomorphism in Corollary \ref{Cor:canonical homomorphism for the filtered case}, we obtain
\beqn
(\Phi_0^F)_*:\HF_\mathrm{loc}(\AA^\HH)\cong\HF^{(-a,0)}(\AA^\HH)\pf\HF^{(-a+||F||_-,|F||_-)}(\AA^\HH_F)\cong\HF_\mathrm{loc}(\AA^\HH_F).
\eeq
On the other hand, we also have
\beqn
(\Phi_F^0)_*:\HF^{(-a+||F||_-,|F||_-)}(\AA^\HH_F)\pf\HF^{(-a+||F||,||F||)}(\AA^\HH)\cong\RFH_\mathrm{loc}(\Sigma,M).
\eeq
Using a homotopy of homotopies $D_s^r(t,x)=\varphi_r(s)F(t,x)$, we deduce
\beqn
(\Phi_F^0)_*\circ(\Phi_0^F)_*=\id_{\HF_\mathrm{loc}(\AA^\HH)}.
\eeq
Therefore $(\Phi_0^F)_*$ is injective and the proposition follows with $\iota=(\Phi_0^F)_*\circ\Theta$.
\end{proof}
\noindent\textbf{Proof of Theorem C.} It directly follows from Proposition \ref{Prop:betti bumber inequalities} and Proposition \ref{Prop:injective homomorphism}.
\hfill $\square$

\end{document}